\documentclass[a4paper,11pt]{article}%
\usepackage{amsmath}
\usepackage{amsfonts}
\usepackage{amssymb}
\usepackage{graphicx}%
\newtheorem{theorem}{Theorem}

\newtheorem{corollary}[theorem]{Corollary}

\newtheorem{definition}[theorem]{Definition}
\newtheorem{propdef}[theorem]{Proposition-definition}
\newtheorem{theoremdef}[theorem]{Theorem-definition}
\newtheorem{example}[theorem]{Example}

\newtheorem{lemma}[theorem]{Lemma}
\newtheorem{notation}[theorem]{Notation}

\newtheorem{proposition}[theorem]{Proposition}
\newtheorem{remark}[theorem]{Remark}

\begin{document}

\title{{Generalized Integral Operators and Applications}}
\author{S.Bernard \thanks{Universit\'{e} des Antilles et de la Guyane, Laboratoire
AOC, Campus de Fouillole, 97159 Pointe-\`{a}-Pitre, Guadeloupe, E-mail:
severine.bernard@univ-ag.fr}, J.-F.Colombeau \thanks{Institut Fourier,
Universit\'{e} J. Fourier, 100 rue des Maths, BP 74, 38402 St Martin d'Heres,
France, E-mail: jf.colombeau@wanadoo.fr}, A.Delcroix \thanks{Same address as
the first author, E-mail: antoine.delcroix@univ-ag.fr}}
\date{}
\maketitle

\begin{abstract}
We extend the theory of distributional kernel operators to a framework of
generalized functions, in which they are replaced by integral kernel
operators. Moreover, in contrast to the distributional case, we show that
these generalized integral operators can be composed unrestrictedly. This
leads to the definition of the exponential, and more generally entire
functions, of a subclass of such operators.

\end{abstract}


\section{Introduction}

The theory of nonlinear generalized functions
\cite{biagioni,colombeau84,colombeau85,colombeau92,farkas,grosser,nedeljkov},
which appears as a natural extension of the theory of
distributions, seems to be a suitable framework to overcome the
limitations of the classical theory of unbounded operators.

Following a first approach done by D. Scarpalezos in
\cite{scarpalezos}, we introduce a natural concept of generalized
integral kernel operators in this setting. In addition, we show
that these operators are characterized by their kernel. Our
approach has some relationship with the one of
\cite{garetto2,garetto1} but is less restrictive and uses other
technics of proofs. Let us quote that classical operators with
smooth or distributional kernel are represented by generalized
integral kernel operators in the spaces of generalized functions,
through the sheaf embeddings of $\mathcal{C}^{\infty}$ or $\mathcal{D}%
^{\prime}$ into $\mathcal{G}$, the sheaf of spaces of generalized functions.
This shows that our theory is a natural extension of the classical one.

Contrary to the classical case \cite{vokhac}, we show that such
operators can be composed unrestrictedly. This is done for
generalized operators with kernel properly supported belonging to
the classical space of generalized functions ${\mathcal{G}}$ and
for operators\ with kernel in a less usual space
${\mathcal{G}}_{L^{2}}$, constructed from the algebra ${\mathcal{D}}_{L^{2}%
}=H^{\infty}$. This allows to consider their iterate composition and the
question of summation of series of such operators naturally arises. In view of
applications to theoretical physics, this question has been solved for the
exponential, with additional assumptions on the growth of the kernel with
respect to the scaling parameter.\ Two cases have been considered: The case of
operators with compactly supported kernel for which the results have been
announced and partially proved in \cite{BCD2}; The case of operators with
kernels in the above mentioned space ${\mathcal{G}}_{L^{2}}$ for which we give
an application to symmetrical operators.

\section{The mathematical framework}

In order to render the paper almost self contained, we recall some
elements of the theory of generalized numbers and functions
without any proofs. We refer the reader to
\cite{biagioni,colombeau84,colombeau85,colombeau92,farkas,grosser,marti98,marti03,nedeljkov}
for more details (except for subsection \ref{prelim}).

\subsection{The sheaf of algebras of generalized functions}

Let $E$ be a sheaf of topological $\mathbb{K}$-algebras on a topological space
$X$ ($\mathbb{K}=\mathbb{R}$ or $\mathbb{C}$). As in \cite{marti98}, we assume
that $E$ satisfies the two following properties:

\begin{itemize}
\item[(i)] For each open subset $\Omega$ of $X$, the algebra $E(\Omega)$ is endowed
with a family of semi-norms $\mathcal{P}$$(\Omega)=(p_{i})_{i\in I(\Omega)}$,
which gives to $E(\Omega)$ a structure of topological vector space and
satisfies
\begin{equation}
\forall i\in I(\Omega),~\exists(j,k,C)\in I(\Omega)^{2}\times\mathbb{R}%
_{+}^{\ast}~/~\forall f,g\in E(\Omega)~,~p_{i}(fg)\leq
Cp_{j}(f)p_{k}(g);
\label{GIO-ContP}%
\end{equation}

\item[(ii)] For two open subsets $\Omega_{1}$ and $\Omega_{2}$ of $X$ such that
$\Omega_{1}\subset\Omega_{2}$, one has
$$\forall i\in I(\Omega_{1})~,~\exists j\in I(\Omega_2)~/~\forall u\in E(\Omega_{2})~,~p_{i}%
(u\mid_{~\Omega_{1}})\le p_{j}(u);$$

\item[(iii)] Let
${\mathcal{F}}=(\Omega_{\lambda})_{\lambda\in\Lambda}$ be any
family of open subsets of $X$ with
$\displaystyle\Omega=\cup_{\lambda\in\Lambda}\Omega_{\lambda}$.
Then, for each $p_i\in{\mathcal{P}}(\Omega)$, $i\in I(\Omega)$,
there exists a finite subfamily of ${\mathcal{F}}$:
$\Omega_1,\Omega_2,...,\Omega_{s(i)}$ and corresponding semi-norms
$p_1\in{\mathcal{P}}(\Omega_1)$,
$p_2\in{\mathcal{P}}(\Omega_2)$,...,
$p_{s(i)}\in{\mathcal{P}}(\Omega_{s(i)})$, such that, for any
$u\in E(\Omega)$,
$$p_i(u)\le \max_{1\le i \le s(i)}\left( p_{i}%
(u\mid_{~\Omega_{i}}) \right).$$
\end{itemize}

Set%
\begin{gather*}
\mathcal{H}(E,\mathcal{P})(\Omega)=\left\{  (u_{\varepsilon})_{\varepsilon}\in
E(\Omega)^{(0,1]}~/~\forall i\in I(\Omega),~\exists~n\in\mathbb{N}%
~:~p_{i}(u_{\varepsilon})=O(\varepsilon^{-n})\ \mathrm{as\ }\varepsilon
\rightarrow0\right\} \\
\mathcal{I}(E,\mathcal{P})(\Omega)=\left\{  (u_{\varepsilon})_{\varepsilon}\in
E(\Omega)^{(0,1]}~/~\forall i\in I(\Omega),~\forall~n\in\mathbb{N}%
~:~p_{i}(u_{\varepsilon})=O(\varepsilon^{n})\ \mathrm{as\ }\varepsilon
\rightarrow0\right\}  .
\end{gather*}

As proved in \cite{marti98}, the functor $\mathcal{H}$$(E,$$\mathcal{P}%
$$):\Omega\mapsto$$\mathcal{H}$$(E,$$\mathcal{P}$$)(\Omega)$ is a sheaf of
subalgebras of the sheaf $E^{(0,1]}$, the functor $\mathcal{I}$$(E,$%
$\mathcal{P}$$):\Omega\mapsto$$\mathcal{I}$$(E,$$\mathcal{P}$$)(\Omega)$ is a
sheaf of ideals of $\mathcal{H}$$(E,$$\mathcal{P}$$)$ and the constant factor
sheaf $\mathcal{H}$$(\mathbb{K},|\cdot|)/$$\mathcal{I}$$(\mathbb{K},|\cdot|)$
is exactly the factor ring $\overline{\mathbb{K}}=A/I_{A}$, with%
\begin{gather*}
A=\{(r_{\varepsilon})_{\varepsilon}\in\mathbb{K}^{(0,1]}~/~\exists
~n\in\mathbb{N}~:~|r_{\varepsilon}|=O(\varepsilon^{-n})\ \mathrm{as\ }%
\varepsilon\rightarrow0\}\\
I_{A}=\{(r_{\varepsilon})_{\varepsilon}\in\mathbb{K}^{(0,1]}~/~\forall
~n\in\mathbb{N}~:~|r_{\varepsilon}|=O(\varepsilon^{n})\ \mathrm{as\ }%
\varepsilon\rightarrow0\}.
\end{gather*}
The sheaf of factor algebras $\mathcal{A}$$\left(  E,\mathcal{P}\right)
\mathcal{=H}$$(E,$$\mathcal{P}$$)/$$\mathcal{I}$$(E,$$\mathcal{P}$$)$, is
called a \textit{sheaf of nonlinear generalized functions}.

\begin{remark}
If $E$ is a sheaf of differential algebras then the same holds for
$\mathcal{A}$$(E,$$\mathcal{P}$$)$.
\end{remark}

In the following paragraphs, we will use this definition in the following
particular cases.

\begin{example}
\label{GIORing}We define $\overline {\mathbb{C}}$ (resp.
$\overline{\mathbb{R}}$) to be the factor ring $A/I_{A}$ of
generalized complex (resp. real) numbers.
\end{example}

\begin{example}
\label{GIOExCol}Take the sheaf $E=$ $\mathcal{C}$$^{\infty}$ on $X=\mathbb{R}%
^{d}$ $(d\in\mathbb{N})$, endowed with its usual topology. This
topology can be described, for $\Omega$ an open subset of
$\mathbb{R}^{d}$, by the family
$\mathcal{P}$$(\Omega)=\{p_{K,l}\,;K\Subset\Omega,~l\in\mathbb{N}\}$,
where the notation $K\Subset\Omega$ means that $K$ is a compact
subset included in
$\Omega$ and%
\[
p_{K,l}(f)=\sup_{x\in K,|\alpha|\leq l}|\partial^{\alpha}%
f(x)|,~\mathrm{for\ all}\ f\in\mathcal{C}^{\infty}(\Omega)\text{ (with
}\partial^{\alpha}=\frac{\partial^{\left\vert \alpha\right\vert }f}{\partial
x^{\alpha}}\text{)}.
\]
$\mathcal{A}$$($$\mathcal{C}$$^{\infty},$$\mathcal{P}$$)\left(  \Omega\right)
$ is the algebra of simplified generalized functions, introduced by the second
author \cite{biagioni,colombeau84,colombeau85,colombeau92,grosser}.
\end{example}

\begin{notation}
We set $\mathcal{E}$$_{M}(\Omega)=\mathcal{H}$$($$\mathcal{C}$$^{\infty}%
,$$\mathcal{P}$$)(\Omega)$, $\mathcal{I}$$(\Omega)=\mathcal{I}$$($%
$\mathcal{C}$$^{\infty},$$\mathcal{P}$$)(\Omega)$ and $\mathcal{G}$%
$(\Omega)=\mathcal{A}$$($$\mathcal{C}$$^{\infty},$$\mathcal{P}$$)(\Omega)$. We
shall also write $P_{K}$ instead of $P_{K,0}$ for every compact subset $K$ of
$\Omega$.
\end{notation}

Since $\mathcal{G}$ is a sheaf, the support of a section $u\in$ $\mathcal{G}%
(\Omega)$ is well defined. Let us recall that, for
$\Omega^{\prime}$ an open subset of $\Omega$ and
$u\in\mathcal{G}$$(\Omega)$, the restriction of $u$ to
$\Omega^{\prime}$ is the class in $\mathcal{G}$$(\Omega^{\prime})$
of $\left( u_{\varepsilon\left\vert \Omega^{\prime}\right.
}\right)  _{\varepsilon}$ where $\left(  u_{\varepsilon}\right)
_{\varepsilon}$ is any representative of $u$. We say that $u$ is
null on $\Omega^{\prime}$ if its restriction to $\Omega^{\prime}$
is null in $\mathcal{G}$$(\Omega^{\prime})$. The \emph{support} of
a generalized function $u\in$ $\mathcal{G}$$(\Omega)$ is the
complement in $\Omega$ of the largest open subset of $\Omega$
where $u$ is null.

\begin{notation}
For $\Omega$ an open subset of $\mathbb{R}^{d}$, we will denote by ${\mathcal{G}%
}_{C}(\Omega)$ the set of generalized functions of ${\mathcal{G}}(\Omega)$
with compact support.
\end{notation}

\begin{remark}
Every $f\in\mathcal{G}_{C}\left(  \Omega\right)  $ has a
representative $\left(  f_{\varepsilon}\right)
_{\varepsilon}\in\mathcal{E}_{M}\left( \Omega\right)  $, such that
each $f_{\varepsilon}$ is supported in the same compact set. We
say that such a representative has a \emph{global compact
support}.
\end{remark}

The two following examples will be used in section
\ref{exponentiel} for the definition of the exponential of some
generalized integral operator. In them, we apply the Colombeau
construction to presheaves of algebras in example \ref{GIOGL2}
(resp. vector spaces in example \ref{GIOGL2l}). In these cases,
property (iii) may not be satisfied but the general construction
is still valid, giving presheaves of generalized algebras (resp.
of vector spaces).

\begin{example}
\label{GIOGL2}For $\Omega$ an open subset of $\mathbb{R}^{d}$, we
consider $E(\Omega)=H^{\infty}(\Omega)$ with its usual topology,
defined by the family
$\mathcal{P}$$(\Omega)=\{\Vert\cdot\Vert_{m}^{\prime}\,;~m\geq0\}$,
with
\[
\Vert f\Vert_{m}^{\prime}=\Vert f\Vert_{H^{m}(\Omega)}=\sum_{|\alpha|\leq
m}\Vert\partial^{\alpha}f\Vert_{L^{2}(\Omega)},~for\ all\ f\in H^{\infty
}(\Omega).
\]

\end{example}

\begin{notation}
\label{GIOGL2Not}We set $\mathcal{E}$$_{L^{2}}(\Omega)=\mathcal{H}$%
$(H^{\infty},$$\mathcal{P}$$)(\Omega)$, $\mathcal{I}$$_{L^{2}}(\Omega
)=\mathcal{I}$$(H^{\infty},\mathcal{P}),$$\mathcal{P}$$)(\Omega)$ and
$\mathcal{G}$$_{L^{2}}(\Omega)=\mathcal{A}$$(H^{\infty},$$\mathcal{P}%
$$)(\Omega).$
\end{notation}

When $E\left(  \Omega\right)  $ is only a topological vector space on
$\mathbb{K}$ (that is (\ref{GIO-ContP}) is not necessarly satisfied),
$\mathcal{G}$$(\Omega)$ is defined analogously and is still a module on
$\overline{\mathbb{K}}$.

\begin{example}
\label{GIOGL2l}For $\Omega$ an open subset of $\mathbb{R}^{d}$, we
consider
$E^{i}(\Omega)=L^{i}(\Omega)\cap\mathcal{C}$$^{\infty}(\Omega)$
($i=1$ and $i=2$) with\ the topology given by the norm
$\Vert\cdot\Vert_{i}=\left\Vert \cdot\right\Vert
_{L^{i}(\Omega)}.$
\end{example}

\begin{notation}
We set $\mathcal{E}$$_{l^{i}}(\Omega)=\mathcal{H}$$(L^{i}\cap\mathcal{C}%
^{\infty},\Vert f\Vert_{L^{i}(\Omega)})(\Omega)$, ${\mathcal{I}}_{l^{2}%
}(\Omega)=\mathcal{I}$$(L^{i}\cap\mathcal{C}^{\infty},\Vert f\Vert
_{L^{i}(\Omega)})$ and $\mathcal{G}$$_{l^{i}}(\Omega)=\mathcal{E}$$_{l^{i}%
}(\Omega)/{\mathcal{I}}_{l^{i}}(\Omega)$.
\end{notation}

\subsection{Embeddings of spaces of distributions into spaces of generalized
functions\label{injection}}

Let $\Omega$ be an open subset of $\mathbb{R}^{d}$
($d\in\mathbb{N)}$. The embedding of
${\mathcal{C}}^{\infty}(\Omega)$ into ${\mathcal{G}}(\Omega)$ is
given by the canonical map
\[
\sigma~:~{\mathcal{C}}^{\infty}(\Omega)\rightarrow{\mathcal{G}}(\Omega
)\ \ \ \ \ f\mapsto Cl(f_{\varepsilon})_{\varepsilon}\text{, with
}f_{\varepsilon}=f\text{ for all }\varepsilon\in\left(  0,1\right]  \text{,}%
\]
which is an injective homomorphism of algebras.

An embedding $i_{S}$ of ${\mathcal{D}}^{\prime}(\Omega)$ into ${\mathcal{G}%
}(\Omega)$ such that $i_{S}\mid_{{\mathcal{C}}^{\infty}(\Omega)}=\sigma$ can
be constructed by the two following methods. For the first one \cite{grosser},
one starts from a net $\left(  \rho_{\varepsilon}\right)  _{\varepsilon}$
defined by $\rho_{\varepsilon}\left(  \cdot\right)  =\varepsilon^{-d}%
\rho\left(  \cdot/\varepsilon\right)  $, where $\rho\in{\mathcal{S}%
}(\mathbb{R}^{d})$ satisfies%
\[%
{\textstyle\int}
\rho(x)\,\mathrm{d}x=1\ ;\ \ \forall
m\in\mathbb{N}^{d}\setminus\{0\}\ \ \
{\textstyle\int}
x^{m}\rho(x)\,\mathrm{d}x=0.
\]
An embedding $i_{0}$ of ${\mathcal{E}}^{\prime}(\mathbb{R}^{d})$ in
${\mathcal{G}}(\mathbb{R}^{d})$ is defined by
\[
i_{0}~:~{\mathcal{E}}^{\prime}(\Omega)\rightarrow{\mathcal{G}}(\Omega
)\ \ \ \ T\mapsto
Cl((T\ast\rho_{\varepsilon})_{\mid\Omega})_{\varepsilon}.
\]
From this, for every open subset $\Omega\subset\mathbb{R}^{d}$, an
open covering $(\Omega_{\lambda})_{\lambda}$ of $\Omega$ with
relatively compact open subsets is considered, and
${\mathcal{D}}^{\prime}(\Omega_{\lambda})$ is embedded into
${\mathcal{G}}(\Omega_{\lambda})$ with the help of cutoff
functions and $i_{0}$. Using a partition of unity subordinate to
$(\Omega_{\lambda })_{\lambda}$, the embedding $i_{S}$ of
${\mathcal{D}}^{\prime}(\Omega)$ into ${\mathcal{G}}(\Omega)$ is
constructed by gluing the bits obtained before together. Finally,
it is shown that the embedding $i_{S}$ does not depend on the
choice of $(\Omega_{\lambda})_{\lambda}$ and other material of the
construction, excepted the net
$(\rho_{\varepsilon})_{\varepsilon}$. The second method
\cite{nedeljkov} starts from the same $\left(  \rho
_{\varepsilon}\right)  _{\varepsilon}$ which is slightly modified
by ad hoc cutoff functions. Consider
$\chi\in{\mathcal{D}}(\mathbb{R})$ even such that
\[
0\leq\chi\leq1\,,~\ \ \ \chi\equiv1\ \mathrm{on}\ \bar{B}(0,1)\,,~\ \ \ \chi
\equiv0\ \mathrm{on}\ \mathbb{R}^{d}\setminus B(0,2)
\]
and set
\[
\forall x\in\mathbb{R}^{d},\ \ \forall\varepsilon\in\left(  0,1\right]
,\ \ \Theta_{\varepsilon}(x)=\rho_{\varepsilon}(x)\,\chi(|\ln\varepsilon|x).
\]
One shows that%
\begin{equation}
\left(
{\textstyle\int}
\Theta_{\varepsilon}(x)\,\mathrm{d}x-1\right)  _{\varepsilon}\in{\mathcal{I}%
}(\mathbb{R})\ ;\ \ \forall m\in\mathbb{N}^{d}\setminus\{0\},\ \ \left(
{\textstyle\int}
x^{m}\Theta_{\varepsilon}(x)\,\mathrm{d}x\right)  _{\varepsilon}%
\in{\mathcal{I}}(\mathbb{R}). \label{ScarpMol}%
\end{equation}
Set $\Gamma_{\varepsilon}=\left\{  x\in\Omega~/~d(x,\mathbb{R}^{d}%
\setminus\Omega)\geq\varepsilon~,~d(x,0)\leq1/\varepsilon\right\}  $ and
consider $(\gamma_{\varepsilon})_{\varepsilon}\in{\mathcal{D}}(\mathbb{R}%
^{d})^{(0,1]}$ such that
\[
\forall\varepsilon\in(0,1]\ ,\ \ 0\leq\gamma_{\varepsilon}\leq1,\ \ \gamma
_{\varepsilon}\equiv1\ \mathrm{on\ }\Gamma_{\varepsilon}.
\]
Then the map%
\[
\mathcal{D}^{\prime}(\Omega)\rightarrow{\mathcal{G}}(\Omega)\ \ \ \ T\mapsto
Cl(\gamma_{\varepsilon}T\ast\Theta_{\varepsilon})_{\varepsilon}%
\]
is equal to $i_{S}$ \cite{delcroix1}. (This last proof uses mainly
(\ref{ScarpMol}); The additional cutoff $\left(
\gamma_{\varepsilon}\right) _{\varepsilon}$, which is such that
$\gamma_{\varepsilon}T\mapsto T$ in
${\mathcal{D}}^{\prime}(\Omega)$ as $\varepsilon\rightarrow0$, is
needed to obtain a well defined net
$(\gamma_{\varepsilon}T\ast\Theta_{\varepsilon })_{\varepsilon}$.)

\subsection{Integration of generalized functions\label{intfg}}

We shall use integration of generalized functions on compact sets or
integration of generalized functions having compact support.

Let $K$ be a given compact subset of $\Omega$ and $u$ an element of
$\mathcal{G}$$(\Omega)$. The integral of $u$ on $K$, denoted by $\int
_{K}u(x)\,\mathrm{d}x$, is the class, in $\overline{\mathbb{C}}$ of the
integral on $K$ of any representative of $u$. (This class does not depend on
the choice of the representative of $u$.)

The integral of a generalized function having a compact support is an
immediate extension of the previous case. Indeed, if $u\in$ $\mathcal{G}%
$$(\Omega)$ has a compact support $K$, let $K_{1}\subset K_{2}$ be two compact
subsets of $\Omega$ such that $K$ is contained in the interior of $K_{1}$.
Then, it can be shown that
\[
\int_{K_{2}\setminus\overset{\circ}{K_{1}}}u(x)\,\mathrm{d}x=0\ \mathrm{in}%
\ \overline{\mathbb{C}}.
\]
Therefore, $\int_{K_{1}}u(x)\,\mathrm{d}x=\int_{K_{2}}u(x)\,\mathrm{d}x$ and
this value is denoted by $\int_{\Omega}u(x)\,\mathrm{d}x$. ($\Omega$ is
omitted in the sequel if no confusion may arise.)\smallskip

We shall also consider integration on the space $\mathcal{G}$$_{l^{1}}%
(\Omega)$. The integral of $u\in\mathcal{G}$$_{l^{1}}(\Omega)$, denoted by
$\int_{\Omega}u(x)\,\mathrm{d}x$, is the class, in $\overline{\mathbb{C}}$ of
the integral on $\Omega$ of any representative of $u$. (This class does not
depend on the choice of the representative of $u$.) \smallskip

It follows immediately from the definitions that the integral of a
generalized function on a set of measure zero is equal to zero,
the integral of a null generalized function is equal to zero and
that the classical formulas of integration by parts, change of
variables, change in order of integration (Fubini's theorem),
\ldots\ are valid for the integration of generalized functions.

\subsection{Generalized parameter integrals\label{prelim}}

Let $X$ (resp. $Y$) be an open subset of $\mathbb{R}^{m}$ (resp.
$\mathbb{R}^{n}$). We denote by $\mathcal{G}$$_{ps}(X\times Y)$
the set of generalized functions $g$ of $\mathcal{G}$$(X\times Y)$
\emph{properly supported} in the following sense:
\begin{equation}
\forall~O_{1}\ \mathrm{relatively\ compact\ open\ subset\ of\ }X,~\exists
K_{2}\Subset Y~/~supp~g\cap(O_{1}\times Y)~\subset O_{1}\times K_{2}.
\label{hyp}%
\end{equation}
Clearly, $\mathcal{G}$$_{ps}(X\times Y)$ is a subalgebra of $\mathcal{G}%
$$(X\times Y)$.

\begin{lemma}
\label{GIOLmLoc}Let $g$ be in $\mathcal{G}$$_{ps}(X\times Y)$. For $V$
relatively compact open subset of $X$, there exists $W$ relatively compact
open subset of $Y$ such that $supp~g\cap(V\times Y)~\subset V\times
W$.\newline For all $\varepsilon\in(0,1]$ and $x\in V$, we set $G_{\varepsilon
}(x)=\int_{W}g_{\varepsilon}(x,y)\,\mathrm{d}y$, where $(g_{\varepsilon
})_{\varepsilon}$ is a representative of $g$. The net $(G_{\varepsilon
})_{\varepsilon}$ belongs to $\mathcal{E}$$_{M}(V)$ and its class, denoted by
$G$, is an element of $\mathcal{G}$$(V)$ which does not depend on the choice
of the representative of $g$ and of $W$.
\end{lemma}

\underline{{\bf{Proof}}}. First, the existence of $W$ is due to
the hypothesis (\ref{hyp}). Then, for all $\varepsilon\in(0,1]$,
$G_{\varepsilon}$ is well defined and of class
$\mathcal{C}$$^{\infty}$ by the usual regularity theorems. (Note
that $G_{\varepsilon}$ is the integral of a
$\mathcal{C}$$^{\infty}$-function on a relatively compact open
subset.)

We first show that $(G_{\varepsilon})_{\varepsilon}$ belongs to $\mathcal{E}%
$$_{M}(V)$. Let $K_{1}$ be a compact subset of $V$ and $\alpha$ be in
$\mathbb{N}^{m}$. For $x\in K_{1}$, one has%
\begin{align*}
\left.  \left\vert \partial^{\alpha}G_{\varepsilon}(x)\right\vert
\leq\left\vert \int_{W}\partial^{\alpha}g_{\varepsilon}(x,y)\,\mathrm{d}%
y\right\vert \right.   &  \leq Vol(W)\sup_{x\in K_{1},y\in\bar{W}}%
|\partial^{\alpha}g_{\varepsilon}(x,y)|\\
&  \leq Vol(W)\,C\,\varepsilon^{-q},\mbox{ for }\varepsilon\mbox{
small enough},
\end{align*}
for some $C>0$ and $q\in\mathbb{N}$, where $Vol(W)$ denotes the volume of $W$.

Let us verify that $G$ does not depend on the choice of the representative of
$g$ and on the one of $W$. According to M. Grosser \textit{et al.}
(\cite{grosser}, theorem 1.2.3), it is enough to consider estimates of order
zero, what we will do in the following. Let $(g_{\varepsilon}^{1}%
)_{\varepsilon}$ and $(g_{\varepsilon}^{2})_{\varepsilon}$ be two
representatives of $g$. As previously, we can define $(G_{\varepsilon}%
^{1})_{\varepsilon}$ and $(G_{\varepsilon}^{2})_{\varepsilon}$ in
$\mathcal{E}$$_{M}(V)$. We have to show that $(G_{\varepsilon}^{1}%
-G_{\varepsilon}^{2})_{\varepsilon}$ is in $\mathcal{I}$$(V)$. Let $K_{1}$ be
a compact subset of $V$. For $x\in K_{1}$, one has%
\[
\left\vert G_{\varepsilon}^{1}(x)-G_{\varepsilon}^{2}(x)\right\vert
=\left\vert \int_{W}(g_{\varepsilon}^{1}-g_{\varepsilon}^{2})(x,y)\,\mathrm{d}%
y\right\vert \leq Vol(W)\sup_{x\in K_{1},y\in\bar{W}}|(g_{\varepsilon}%
^{1}-g_{\varepsilon}^{2})(x,y)|.
\]
As $(g_{\varepsilon}^{1}-g_{\varepsilon}^{2})_{\varepsilon}$ belongs to
$\mathcal{I}$$(X\times Y)$, we get from the previous estimate that $P_{K_{1}%
}\left(  G_{\varepsilon}^{1}-G_{\varepsilon}^{2}\right)  =O\left(
\varepsilon^{n}\right)  $ as $\varepsilon\rightarrow0$ for all $n\in
\mathbb{N}$.

Consider $W_{1}$ and $W_{2}$ relatively compact open subsets of $Y$ such that,
$supp~g\cap(V\times Y)~\subset V\times W_{i}$ for $i=1,2$ with, for example
$W_{1}\subset W_{2}$. For all $\varepsilon\in(0,1]$, $x\in V$ and $i=1,2$, we
set $G_{\varepsilon}^{i}(x)=\int_{W_{i}}g_{\varepsilon}(x,y)\,\mathrm{d}y$
where $(g_{\varepsilon})_{\varepsilon}$ is a representative of $g$. Then
$(G_{\varepsilon}^{1})_{\varepsilon}$ and $(G_{\varepsilon}^{2})_{\varepsilon
}$ belong to $\mathcal{E}$$_{M}(V)$. Let $K_{1}$ be a compact subset of $V$.
For $x\in K_{1}$, one has%
\[
\left\vert (G_{\varepsilon}^{1}-G_{\varepsilon}^{2})(x)\right\vert
=\left\vert \int_{W_2\setminus
W_1}g_{\varepsilon}(x,y)\,\mathrm{d}y\right\vert \leq
Vol(W_2\setminus W_1)\sup_{x\in K_{1},y\in\overline{W_2\setminus
W_1}}|g_{\varepsilon}(x,y)|.
\]
As $supp~g\cap(V\times Y)~\subset V\times W$, the restriction of
$g$ to $V\times(W_2\setminus W_1)$ is null.\ Therefore, the
previous estimate shows that $P_{K_{1}}\left(
G_{\varepsilon}^{1}-G_{\varepsilon}^{2}\right) =O\left(
\varepsilon^{n}\right)  $ as $\varepsilon\rightarrow0$ for all
$n\in\mathbb{N}$.

\begin{lemma}
\label{GIOLmGlob}Let $g$ be in $\mathcal{G}$$_{ps}(X\times Y)$, $(V_{i})_{i\in
I}$ be a family of relatively compact open subsets of $X$ such that
$\cup_{i\in I}V_{i}=X$ and define $G_{i}\in$ $\mathcal{G}$$(V_{i})$ as in
lemma \ref{GIOLmLoc}. Then, there exists $G\in$ $\mathcal{G}$$(X)$ such that
the restriction of $G$ to $V_{i}$ is equal to $G_{i}$, for all $i\in I$.
Moreover, $G$ only depends on $g$, but not on $(V_{i})_{i\in I}$.
\end{lemma}

\underline{{\bf{Proof}}}. For $i\not =j$ such that $V_{i}\cap
V_{j}\neq\varnothing$, we remark that $V_{i}\cup V_{j}$ is a
relatively compact open subset of $X$. There exists $W$ a
relatively compact open subset of $Y$ such that $supp~g\cap(\left(
V_{i}\cup V_{j}\right)  \times Y)~\subset\left(  V_{i}\cup
V_{j}\right) \times W$. According to lemma \ref{GIOLmLoc}, we can
define $\Phi=Cl\left( \Phi_{\varepsilon}\right)
_{\varepsilon}\in\mathcal{G}$$(V_{i}\cup V_{j})$,
with%
\[
\forall x\in V_{i}\cup V_{j},\ \ \Phi_{\varepsilon}\left(  x\right)  =\int
_{W}g_{\varepsilon}(x,y)\,\mathrm{d}y.
\]
Then, $\left(  \Phi_{\varepsilon\left\vert V_{i}\right.  }\right)
_{\varepsilon}$ (resp. $\left(  \Phi_{\varepsilon\left\vert
V_{j}\right. }\right)  _{\varepsilon}$) is a\ representative of
$G_{i}$ (resp. $G_{j}$) since those representatives depend neither
on the representative of $g$ nor on the choice of appropriate $W$.
Then $G_{i\left\vert V_{i}\cap V_{j}\right. }=G_{j\left\vert
V_{i}\cap V_{j}\right.  }$. Thus, $\left(  G_{i}\right) _{i\in I}$
is a coherent family, which implies the existence of $G$ since
$\mathcal{G}$$(X)$ is a sheaf. The proof of the independence of
$G$ with respect to $(V_{i})_{i\in I}$ follows the same
lines.\medskip

Lemmas \ref{GIOLmLoc} and \ref{GIOLmGlob} give immediately the following:

\begin{proposition}
\label{intgps}For $g$ in $\mathcal{G}$$_{ps}(X\times Y)$, there exists $G\in$
$\mathcal{G}$$(X)$ such that, for all relatively compact open subset $O_{1}$
of $X$,%
\[
G_{\left\vert O_{1}\right.  }=Cl\left(  \,(x\mapsto\int_{K_{2}}g_{\varepsilon
}(x,y)\,\mathrm{d}y)_{\left\vert O_{1}\right.  }\right)  _{\varepsilon},
\]
where $\left(  g_{\varepsilon}\right)  $ is a representative of $g$ and
$K_{2}\Subset Y$ is such that $supp\,g\cap(O_{1}\times Y)~\subset O_{1}\times
K_{2}.$
\end{proposition}

\begin{notation}
\label{NotIntPar}By a slight abuse of notation, we shall set
$G=\int_{Y}g(\cdot ,y)\,\mathrm{d}y$ or $G\left(  \cdot_{1}\right)
=\int_{Y}g(\cdot _{1},y)\,\mathrm{d}y$. We shall omit the set on
which the integration is performed when no confusion may arise.
\end{notation}

\begin{example}
Proposition \ref{intgps} can be used to define the Fourier's transform of a
compactly supported generalized function. Indeed, if $u$ is in ${\mathcal{G}%
}_{C}(Y)$ then the generalized function $\{(x,y)\mapsto
e^{-ixy}u(y)\}$ is in
${\mathcal{G}}_{ps}(X\times Y)$, then $\hat{u}(x)=\int e^{-ixy}u(y)\,\mathrm{d}%
y$ is well defined.
\end{example}

\section{Generalized integral operators}

In this section, we introduce the notion of generalized integral
operator and study their basic properties. The reader can find
another approach in \cite{garetto2,garetto1}. Let $X$ (resp. $Y$)
be an open subset of $\mathbb{R}^{m}$ (resp. $\mathbb{R}^{n}$).

\begin{definition}
\label{GIOdefGIO}Let $H$ be in $\mathcal{G}$$_{ps}(X\times Y)$. We call
\emph{generalized integral operator} the map
\[
\widehat{H}~:~%
\begin{array}
[t]{ccl}%
\mathcal{G}(Y) & \rightarrow & \mathcal{G}(X)\\
f & \mapsto & \widehat{H}(f)=\int H(\cdot,y)f(y)\,\mathrm{d}y,
\end{array}
\]
with the meaning introduced in proposition \ref{intgps} and notation
\ref{NotIntPar}. We say that $H$ is the \emph{kernel} of the generalized
integral operator $\widehat{H}$.
\end{definition}

This map is well defined due to proposition \ref{intgps} since the application
$H(\cdot_{1},\cdot_{2})f(\cdot_{2})$ is clearly is in $\mathcal{G}$%
$_{ps}(X\times Y)$.

\begin{remark}
If $H\in$ $\mathcal{G}$$(X\times Y)$ has a compact support then $H$ satisfies
(\ref{hyp}) and $\widehat{H}$ is well defined. Furthermore, the definition of
$\widehat{H}$ does not need to refer to proposition \ref{intgps} in this case.
Indeed, if $H$ is in $\mathcal{G}_{C}(X\times Y)$ with $supp\,H\subset
\mathring{K}_{1}\times\mathring{K}_{2}$ ($K_{1}\Subset X$, $K_{2}\Subset Y$)
and $f$ in $\mathcal{G}(Y)$, we have
\[
\widehat{H}(f)=Cl\left(  x\mapsto\int_{K_{2}}H_{\varepsilon}%
(x,y)f_{\varepsilon}(y)\,\mathrm{d}y\right)  _{\varepsilon}%
\]
where $\left(  H_{\varepsilon}\right)  _{\varepsilon}$ (resp. $\left(
f_{\varepsilon}\right)  _{\varepsilon}$) is any representative of $H$ (resp.
$f$). Furthermore, as $supp\,H\subset K_{1}\times K_{2}$, we have $H(\cdot
_{1},\cdot_{2})f(\cdot_{2})_{\left\vert \left(  X\backslash K_{1}\right)
\times Y\right.  }=0$, hence $\widehat{H}(f)_{\left\vert X\backslash
K_{1}\right.  }=0$ and $supp\,\widehat{H}(f)\subset K_{1}$.\ (The proof uses
arguments similar to the one of lemma \ref{GIOLmLoc}.) Finally, the image of
$\widehat{H}$ is included in $\mathcal{G}_{C}$$(X)$ and, more precisely, in
$\left\{  g\in\mathcal{G}_{C}(X)~/~supp~g\subset K_{1}\right\}  .$
\end{remark}

\begin{remark}
If $H$ is in $\mathcal{G}$$(X\times Y)$ without any other hypothesis, we can
define a map $\widehat{H}$~: $\mathcal{G}$$_{C}(Y)\rightarrow$ $\mathcal{G}%
$$(X)$ in the same way.\ Indeed, for all $f$ in $\mathcal{G}$$_{C}(Y)$ with
$supp~f=K_{2}$ and for all $O_{1}$ relatively compact open subset of $X$,
$supp~H(\cdot_{1},\cdot_{2})f(\cdot_{2})\cap(~O_{1}\times Y)~\subset
O_{1}\times K_{2}$, that is $H(\cdot_{1},\cdot_{2})f(\cdot_{2})$ is in
$\mathcal{G}$$_{ps}(X\times Y)$. In this case, the generalized integral
operator can be defined globally since $f$ has a representative with global
compact support.
\end{remark}

This remark leads us to make the link between the classical theory of integral
operators acting on $\mathcal{D}(Y)$ and the generalized one. This is detailed
in section \ref{SecClassTh}.

\begin{remark}
\label{SKT}In all previous cases, $\widehat{H}$ is a linear map of
$\overline{\mathbb{C}}$-modules. This holds also for the map\ \ $\widehat
{}~:~\mathcal{G}$$_{ps}(X\times Y)\rightarrow$$\mathcal{L}$$($$\mathcal{G}%
$$(Y),$$\mathcal{G}$$(X))$, which associates $\widehat{H}$ to $H$. Moreover,
$\widehat{H}$ is continuous for the sharp topologies
\cite{scarpalezos}. Conversely, the third author showed in
\cite{delcroix2} that any continuous linear map from
${\mathcal{G}}_{C}(Y)$ to ${\mathcal{G}}(X)$, satisfying
appropriate growth hypotheses with respect to the regularizing
parameter $\varepsilon$, can be written as a generalized integral
operator, giving a Schwartz kernel type theorem in the framework
of integral generalized operators.
\end{remark}

\begin{example}
The identity map of subspaces of compactly generalized functions
with limited growth \cite{delcroix2} admits as kernel
\[
\Phi=Cl\left(  (x,y)\mapsto\Theta_{\varepsilon}(x-y)\right)  _{\varepsilon},
\]
where $(\Theta_{\varepsilon})_{\varepsilon}$ is defined in paragraph
\ref{injection}.
\end{example}

\begin{theorem}
\label{caracker}%
$\mathbf{{(\mbox{Characterization of generalized integral operators by their kernel})}%
}$ One has $\widehat{H}=0$ if and only if $H=0$.
\end{theorem}

A first proof, due to the second author, is based on embeddings of Sobolev's
spaces in spaces of smooth functions. The proof given below is due to V.
Valmorin (personal communication) and uses the following:

\begin{lemma}
Let $\Omega$ be an open subset of $\mathbb{R}^{d}$ and $K$ be a
compact of $\Omega$, of diameter $D$. If $\Phi$ belongs to
$\mathcal{D}$$_{K}(\Omega)$ then
\[
\sup_{y\in K}|\Phi(y)|\leq\left(  D^{d}\int_{\Omega}|\frac{\partial^{d}%
}{\partial y_{1}\partial y_{2}...\partial y_{d}}\Phi(y)|^{2}%
\,\mathrm{\,\mathrm{d}}y\right)  ^{1/2}.
\]

\end{lemma}

\underline{{\bf{Proof}}}. First, let us assume that $H=0$ in
$\mathcal{G}$$(X\times Y)$ and fix $f\in$ $\mathcal{G}$$(Y)$. For
any $O_{1}$ relatively compact open subset of $X$, there exists
$K_{2}\Subset Y$ such that $supp\,g\cap(O_{1}\times Y)~\subset
O_{1}\times K_{2}$ and
\[
\widehat{H}\left(  f\right)  _{\left\vert O_{1}\right.  }=Cl(\,(x\mapsto
\int_{K_{2}}H_{\varepsilon}(x,y)f_{\varepsilon}(y)\,\mathrm{\,\mathrm{d}%
}y)_{\left\vert O_{1}\right.  }\,)_{\varepsilon}\,,
\]
where $\left(  H_{\varepsilon}\right)  _{\varepsilon}$ (resp. $\left(
f_{\varepsilon}\right)  _{\varepsilon}$) is any representative of $H$ (resp.
$f$). As $H$ is null, the map $H_{\varepsilon}(\cdot_{1},\cdot_{2})f(\cdot
_{2})$ is null on $O_{1}\times Y$.\ Thus $\widehat{H}\left(  f\right)
_{\left\vert O_{1}\right.  }$ is null and, by sheaf properties, $\widehat
{H}\left(  f\right)  $ is null.

Conversely, suppose that $\widehat{H}=0$. In order to prove that
$H=0$ in $\mathcal{G}$$(X\times Y)$, we shall prove that
$H_{\left\vert O_{1}\times Y\right.  }=0$ in
$\mathcal{G}$$(O_{1}\times Y)$, for any $O_{1}$ relatively compact
open subset of $X$ and conclude by using the sheaf properties of
$\mathcal{G}$$(\cdot)$. Let $K_{1}$ and $K_{2}$ two compacts
subsets of $O_{1}$ and $Y$ respectively. From (\ref{hyp}), we can
find $W$ a relatively compact open subset of $Y$ such that
$K_{2}\subset W$ and $supp~H\cap(O_{1}\times Y)~\subset
O_{1}\times W$. Let $(H_{\varepsilon})_{\varepsilon}$ be a
representative of $H$ and set
$\varphi_{\varepsilon,x}(y)=H_{\varepsilon }(x,y)\rho(y)$, for all
$y\in Y$ and $x\in O_{1}$, where $\rho$ is a
$\mathcal{C}$$^{\infty}$-function on $Y$ such that $\rho=1$ on $W$
and $supp~\rho\subset O_{2}$, with $O_{2}$ a relatively compact
open subset of
$Y$. Thus, for all $x$ in $O_{1}$, $\varphi_{\varepsilon,x}\in$ $\mathcal{D}%
$$_{\bar{O_{2}}}(Y)$. This implies
\[
\sup_{y\in\bar{O_{2}}}|\varphi_{\varepsilon,x}(y)|\leq\left(  diam(\bar{O_{2}%
})^{n}\int_{O_{2}}|\frac{\partial^{n}}{\partial y^{n}}\varphi_{\varepsilon
,x}(y)|^{2}\,\mathrm{d}y\right)  ^{1/2},
\]
where $\frac{\partial^{n}}{\partial y^{n}}$ is the derivative $\frac
{\partial^{n}}{\partial y_{1}\partial y_{2}...\partial y_{n}}$ and
$diam(\cdot)$ denotes the diameter. As $K_{2}\subset W\subset\bar{O_{2}}$, we
have
\[
\sup_{y\in K_{2}}|\varphi_{\varepsilon,x}(y)|\leq\left(  diam(\bar{O_{2}}%
)^{n}\int_{O_{2}}|\frac{\partial^{n}}{\partial y^{n}}\varphi_{\varepsilon
,x}(y)|^{2}\,\mathrm{d}y.\right)  ^{1/2}.
\]
Set
\[
\psi_{\varepsilon}(x)=\int_{O_{2}}|\frac{\partial^{n}}{\partial y^{n}}%
\varphi_{\varepsilon,x}(y)|^{2}\,\mathrm{d}y~,~\forall x\in O_1.
\]
Since $H_{\varepsilon}$ and $\rho$ are $\mathcal{C}$$^{\infty}$-functions,
$\psi_{\varepsilon}$ is continuous on $K_{1}$. Therefore, $\psi_{\varepsilon}$
has its maximum at a point $x(\varepsilon)\in K_{1}$. Consequently,
\[
\sup_{x\in K_{1},y\in K_{2}}|H_{\varepsilon}(x,y)|=\sup_{x\in K_{1},y\in
K_{2}}|\varphi_{\varepsilon,x}(y)|\leq\sqrt{diam(\bar{O_{2}})^{n}%
\psi_{\varepsilon}(x(\varepsilon))}.
\]
By choosing $f_{\varepsilon}=\overline{\frac{\partial^{n}}{\partial y^{n}%
}\varphi_{\varepsilon,x(\varepsilon)}(\cdot)}$, one has $(f_{\varepsilon
})_{\varepsilon}$ in $\mathcal{E}$$_{M}(Y)$, since $x(\varepsilon)\in K_{1}$,
so its class $f$ is an element of $\mathcal{G}$$(Y)$, as well as
$\frac{\partial^{n}}{\partial y^{n}}f$. Since $W\subset O_{2}$, one has
\begin{align*}
(\widehat{H}(\frac{\partial^{n}}{\partial y^{n}}f))_{\varepsilon}%
(x(\varepsilon))  &  =\int_{O_{2}}H_{\varepsilon}(x(\varepsilon
),y)\frac{\partial^{n}}{\partial y^{n}}f_{\varepsilon}(y)\,\mathrm{d}y\\
&  =\int_{O_{2}}\varphi_{\varepsilon,x(\varepsilon)}(y)\frac{\partial^{n}%
}{\partial y^{n}}\overline{\frac{\partial^{n}}{\partial y^{n}}\varphi
_{\varepsilon,x(\varepsilon)}(y)}\,\mathrm{d}y\\
&  =(-1)^{n}\int_{O_{2}}|\frac{\partial^{n}}{\partial y^{n}}\varphi
_{\varepsilon,x(\varepsilon)}(y)|^{2}\,\mathrm{d}y=(-1)^{n}\psi_{\varepsilon
}(x(\varepsilon)),
\end{align*}
since $(x,y)\mapsto\varphi_{\varepsilon,x}(y)$ is another
representative of $H$ on $O_1\times Y$. As
$\frac{\partial^{n}f}{\partial y^{n}}$ belongs to
$\mathcal{G}$$(Y)$, we have $\widehat{H}(\frac{\partial^{n}f}{\partial y^{n}%
})=0$ in $\mathcal{G}$$(X)$, thus its representative is in $\mathcal{I}%
$$(O_{1})$. Since $(x(\varepsilon))_{\varepsilon}$ is bounded in $K_{1}$,
$(\widehat{H}(\frac{\partial^{n}f}{\partial y^{n}}))_{\varepsilon}%
(x(\varepsilon))$ is in $I_{A}$ and the previous equality implies
that $(\psi_{\varepsilon}(x(\varepsilon)))_{\varepsilon}\in
I_{A}$, that is
$\psi_{\varepsilon}(x(\varepsilon))=O(\varepsilon^{n})$ as
$\varepsilon \rightarrow0$, with $n\in\mathbb{N}$. Furthermore
\[
\sup_{(x,y)\in K}|H_{\varepsilon}(x,y)|\leq\sup_{(x,y)\in K_{1}\times K_{2}%
}|H_{\varepsilon}(x,y)|\leq\sqrt{diam(\bar{O_{2}})^{n}\psi_{\varepsilon
}(x(\varepsilon))}.
\]
Hence $\sup_{(x,y)\in K}|H_{\varepsilon}(x,y)|=O(\varepsilon^{p})$
as $\varepsilon\rightarrow0$, for all $p\in\mathbb{N}$. Thus
$(H_{\varepsilon })_{\varepsilon}$ is in
$\mathcal{I}$$(O_{1}\times Y),$ which ends the proof.

\begin{corollary}
The linear map $\widehat{}$ , defined in remark \ref{SKT}, is injective.
\end{corollary}

\begin{remark}
\label{RemOpsGL2}If $H$ is in ${\mathcal{G}}_{L^{2}}(X\times Y)$, then $H$ can
be embedded in ${\mathcal{G}}(X\times Y)$ using Sobolev's embeddings, but $H$
may not be properly supported. Nevertheless, we can define a generalized
integral operator acting on ${\mathcal{G}}_{L^{2}}$ type spaces as follows:
\[
\widehat{H}:%
\begin{array}
[t]{ccl}%
{\mathcal{G}}_{L^{2}}(Y) & \rightarrow & {\mathcal{G}}_{L^{2}}(X)\\
f & \mapsto & Cl\left(  x\mapsto\int_{Y}H_{\varepsilon}(x,y)f_{\varepsilon
}(y)\,\mathrm{d}y\right)  _{\varepsilon},
\end{array}
\]
independently of the representative $(H_{\varepsilon})_{\varepsilon}$ (resp.
$(f_{\varepsilon})_{\varepsilon}$) of $H$ (resp. $f$).
\end{remark}

Indeed, set $\Phi_{\varepsilon}\left(  x\right)  =\int_{Y}H_{\varepsilon
}(x,y)f_{\varepsilon}(y)\,\mathrm{d}y$ for $\varepsilon$ in $\left(
0,1\right]  $ and $x\in X$. The function $\Phi_{\varepsilon}$ is well defined,
since $\left(  H_{\varepsilon}\right)  _{\varepsilon}\in\mathcal{E}_{L^{2}%
}(X\times Y)$ and $\left(  f_{\varepsilon}\right)  _{\varepsilon}%
\in\mathcal{E}$$_{L^{2}}(Y)$. Moreover, we have, for all $x\in X$ and
$\varepsilon$ in $\left(  0,1\right]  $
\[
\left\vert \Phi_{\varepsilon}\left(  x\right)  \right\vert \leq\left\vert
\int_{Y}H_{\varepsilon}(x,y)f_{\varepsilon}(y)\,\mathrm{d}y\right\vert
\leq\left\Vert H_{\varepsilon}(x,\cdot)\right\Vert _{2}\left\Vert
f_{\varepsilon}\right\Vert _{2}.
\]
\ Thus $\Phi_{\varepsilon}$ is in $L^{2}(X)$, with $\left\Vert \Phi
_{\varepsilon}\right\Vert _{2}\leq\left\Vert H_{\varepsilon}\right\Vert
_{2}\left\Vert f_{\varepsilon}\right\Vert _{2}$. As for all $\alpha
\in\mathbb{N}^{n}$, $\partial^{\alpha}\Phi_{\varepsilon}$ exists and satisfies
$\left\Vert \partial^{\alpha}\Phi_{\varepsilon}\right\Vert _{2}\leq\left\Vert
\partial_{x}^{\alpha}H_{\varepsilon}\right\Vert _{2}\left\Vert f_{\varepsilon
}\right\Vert _{2}$: $\left(  \Phi_{\varepsilon}\right)  _{\varepsilon}$ is in
$\mathcal{E}$$_{L^{2}}(X)$. A straightforward computation shows that
$Cl\left(  \Phi_{\varepsilon}\right)  _{\varepsilon}$ does not depend on the
representative of $H$ and $f$.\ Thus, the operator $\widehat{H}$ is well defined.

\section{Link with the classical theory and regularity
properties\label{SecClassTh}}

In this section, we compare our definition to the classical one, that is when
$H$ is a $\mathcal{C}^{\infty}$-function or a distribution. Let $X$ (resp.
$Y$) be an open subset of $\mathbb{R}^{m}$ (resp. $\mathbb{R}^{n}$).

\begin{theorem}
If $h\in{\mathcal{C}}^{\infty}(X\times Y)$ then the diagram
\[%
\begin{array}
[c]{ccc}%
{\mathcal{E}}^{\prime}(Y) & \overset{\widehat{h}}{\rightarrow} & {\mathcal{C}%
}^{\infty}(X)\\
\downarrow i_{S} &  & \downarrow\sigma\\
{\mathcal{G}}_{C}(Y) & \overset{\widehat{\sigma(h)}}{\rightarrow} &
{\mathcal{G}}(X)
\end{array}
\]
is commutative.
\end{theorem}

\underline{{\bf{Proof}}}. We have to prove that
$\sigma\circ\widehat{h}=\widehat{\sigma(h)}\circ i_{S}$. Let $T$
be in $\mathcal{E}$$^{\prime}(Y)$. Due to the local structure of
distributions, there exist $r\in\mathbb{N}$, a finite family
$(f_{\alpha })_{0\leq|\alpha|\leq r}$ $(\alpha\in\mathbb{N}^{n})$
of continuous on $\mathbb{R}^{n}$ having their support contained
in the same arbitrary neighborhood of the support of $T$, such
that $T=\sum_{0\leq|\alpha|\leq r}\partial^{\alpha}f_{\alpha}$. By
the linearity of the operators under consideration, we can assume
that $T=\partial^{\alpha}f$, where $f$ is a continuous function on
$\mathbb{R}^{n}$ whose support is contained in a neighborhood of
the support of $T$. In this proof and the one of the following
theorem, the exponents of the mollifiers are $1$ for the space
$X$, $2$ for $Y$ and none for $X\times Y$. In this case, a
representative of $\sigma\circ\widehat{h}(T)$ is defined, for all
$x\in X$, by
\[
\widehat{h}(T)(x)=\langle T,h(x,\cdot)\rangle=\langle\partial^{\alpha
}f,h(x,\cdot)\rangle=(-1)^{|\alpha|}\langle f,\partial_{y}^{\alpha}%
h(x,\cdot)\rangle=(-1)^{|\alpha|}\int f(y)\partial_{y}^{\alpha}%
h(x,y)\,\mathrm{d}y.
\]
A representative of $\widehat{\sigma(h)}\circ i_{S}(T)$ is $$\int
h(x,y)(T\ast\Theta_{\varepsilon}^{2})(y)\,\mathrm{d}y  =\int
h(x,y)(f\ast\partial^{\alpha}\Theta_{\varepsilon}^{2})(y)\,\mathrm{d}y
=\int\int h(x,y)f(\lambda)\partial^{\alpha}\Theta_{\varepsilon}
^{2}(y-\lambda)\,\mathrm{d}\lambda\mathrm{d}y.$$ As the functions
$f$ and $\partial^{\alpha}\Theta_{\varepsilon}^{2}$ have compact
supports, the two previous integrals are integrals on compacts
sets.\ Thus, we can apply Fubini's theorem and obtain
$$\int\int h(x,y)f(\lambda)\partial^{\alpha}\Theta_{\varepsilon}^{2}
(y-\lambda)d\lambda dy  =(-1)^{|\alpha|}\int\left(
\partial_{y}^{\alpha}h(x,\cdot)\ast
\Theta_{\varepsilon}^{2}\right)
(\lambda)f(\lambda)\,\mathrm{d}\lambda.$$
As the function
$y\mapsto\left(  \partial_{y}^{\alpha}h(x,\cdot)\ast
\Theta_{\varepsilon}^{2}\right)  (y)$ is a representative of $i_{S}%
(\partial_{y}^{\alpha}h(x,\cdot))$ in ${\mathcal{G}}(Y)$ and since
$\partial_{y}^{\alpha}h(x,\cdot)$ is a ${\mathcal{C}}^{\infty}$-function and
$i_{s}\mid_{{\mathcal{C}}^{\infty}}=\sigma$, one has
\[
\left(  \partial_{y}^{\alpha}h(x,\cdot)\ast\Theta_{\varepsilon}^{2}%
-\partial_{y}^{\alpha}h(x,\cdot)\right)  _{\varepsilon}\in{\mathcal{I}}(Y).
\]
Moreover, $f$ is compactly supported so the difference of the
representatives of $\widehat{\sigma(h)}\circ i_{S}(T)$ and
$\sigma\circ\widehat{h}(T)$ is in ${\mathcal{I}}(X)$. Thus
$\widehat{\sigma(h)}\circ i_{S}(T)=\sigma \circ\widehat{h}(T)$ in
${\mathcal{G}}(X)$, which implies the required result.

\begin{definition}
\label{defreg} The kernel $H\in{\mathcal{G}}(X\times Y)$ of a generalized
integral operator is called regular when $\widehat{H}({\mathcal{G}}%
_{C}(Y)\subset{\mathcal{G}}^{\infty}(X)$, where, for all $\Omega$
open subset of $\mathbb{R}^{d}$ ($d\in\mathbb{N}$),
\[
{\mathcal{G}}^{\infty}(\Omega)={\mathcal{E}}^{\infty}(\Omega)/{\mathcal{I}%
}(\Omega).
\]
with
\[
{\mathcal{E}}^{\infty}(\Omega)=\left\{  (u_{\varepsilon})_{\varepsilon}%
\in{\mathcal{C}}^{\infty}(\Omega)^{(0,1]}~/~\forall K\Subset\Omega~,~\exists
n\in\mathbb{N},~\forall l\in\mathbb{N}~,~p_{K,l}(u_{\varepsilon}%
)=O(\varepsilon^{-n})\mbox{, as }\varepsilon\rightarrow0\right\}  .
\]

\end{definition}

The reader can find more details about this algebra in \cite{oberguggenberger}.

\begin{proposition}
If $h$ is in ${\mathcal{C}}^{\infty}(X\times Y)$, then $\sigma(h)$ is regular
in the above sense.
\end{proposition}

\underline{{\bf{Proof}}}. Let $f$ be in ${\mathcal{G}}_{C}(Y)$.
Then there exists $K_{2}$ compact of $Y$ and
$(f_{\varepsilon})_{\varepsilon}$ a representative of $f$ such
that $supp~f_{\varepsilon}\subset K_{2}$. Consequently, a
representative of $\widehat{\sigma(h)}(f)$ is $\left(
\psi_{\varepsilon}\right) _{_{\varepsilon}}$ with
\[
\psi_{\varepsilon}:x\mapsto\int_{K_{2}}h(x,y)f_{\varepsilon}(y)\,\mathrm{d}y
\]
and, for all $K_{1}$ compact of $X$, $\alpha$ in $\mathbb{N}^{m}$, $x$ in
$K_{1}$,
\[
\left\vert \partial_{x}^{\alpha}\left(  \int_{K_{2}}h(x,y)f_{\varepsilon
}(y)\,\mathrm{d}y\right)  \right\vert \leq Vol(K_{2})\sup_{(x,y)\in
K_{1}\times K_{2}}|\partial_{x}^{\alpha}h(x,y)|\sup_{y\in K_{2}}%
|f_{\varepsilon}(y)|\leq C(\alpha)\varepsilon^{-q},
\]
as $\varepsilon$ tends to 0, where $q$ does not depend on
$\alpha$. That shows that $\widehat{\sigma(h)}(f)$ is in
${\mathcal{G}}^{\infty}(X)$.

\begin{theorem}
If $h\in{\mathcal{D}}^{\prime}(X\times Y)$ then the diagram
\[%
\begin{array}
[c]{ccc}%
{\mathcal{D}}(Y) & \overset{\widehat{h}}{\rightarrow} & {\mathcal{D}}^{\prime
}(X)\\
\downarrow\sigma &  & \downarrow i_{S}\\
{\mathcal{G}}_{C}(Y) & \overset{\widehat{i_{S}(h)}}{\rightarrow} &
{\mathcal{G}}(X)
\end{array}
\]
is commutative.
\end{theorem}

\underline{{\bf{Proof}}}. We have to prove that
$i_{S}\circ\widehat{h}=\widehat{i_{S}(h)}\circ\sigma$. Let $f$ be
in $\mathcal{D}$$(Y)$. A representative of $i_{S}\circ\widehat
{h}(f)$ is defined, for all $x\in X$, by
\begin{align*}
\left(  \gamma_{\varepsilon}^{1}\widehat{h}(f)\ast\Theta_{\varepsilon}%
^{1}\right)  (x)  &  =\langle\gamma_{\varepsilon}^{1}\widehat{h}%
(f),\{\xi\mapsto\Theta_{\varepsilon}^{1}(x-\xi)\}\rangle\\
&  =\langle\widehat{h}(f),\{\xi\mapsto\gamma_{\varepsilon}^{1}(\xi
)\Theta_{\varepsilon}^{1}(x-\xi)\}\rangle\\
&  =\langle h,\{\xi\mapsto\gamma_{\varepsilon}^{1}(\xi)\Theta_{\varepsilon
}^{1}(x-\xi)\}\otimes f\rangle.
\end{align*}
A representative of $\widehat{i_{S}(h)}\circ\sigma(f)$ is, for all $x\in X$,%
\begin{multline*}
\int\langle\gamma_{\varepsilon}^{1}\otimes\gamma_{\varepsilon}^{2}%
h,\{(\xi,\eta)\mapsto\Theta_{\varepsilon}^{1}(x-\xi)\Theta_{\varepsilon}%
^{2}(y-\eta)\}\rangle f(y)\,\mathrm{d}y\\
=\langle\gamma_{\varepsilon}^{1}\otimes\gamma_{\varepsilon}^{2}h,\{(\xi
,\eta)\mapsto\Theta_{\varepsilon}^{1}(x-\xi)\int\Theta_{\varepsilon}%
^{2}(y-\eta)f(y)\,\mathrm{d}y\}\rangle\\
=\langle h,\{(\xi,\eta)\mapsto\gamma_{\varepsilon}^{1}(\xi)\Theta
_{\varepsilon}^{1}(x-\xi)\gamma_{\varepsilon}^{2}(\eta)\int\Theta
_{\varepsilon}^{2}(y-\eta)f(y)\,\mathrm{d}y\}\rangle\\
=\langle h,\{(\xi,\eta)\mapsto\gamma_{\varepsilon}^{1}(\xi)\Theta
_{\varepsilon}^{1}(x-\xi)\gamma_{\varepsilon}^{2}(\eta)\left(  \Theta
_{\varepsilon}^{2}\ast f\right)  (\eta)\}\rangle,
\end{multline*}
since $\Theta_{\varepsilon}^{2}$ is even. Consequently, the difference of
these representatives of $i_{S}\circ\widehat{h}(f)$ and $\widehat{i_{S}%
(h)}\circ\sigma(f)$ is equal to
\[
\langle h,\{\xi\mapsto\gamma_{\varepsilon}^{1}(\xi)\Theta_{\varepsilon}%
^{1}(x-\xi)\}\otimes\{f-\gamma_{\varepsilon}^{2}\left(  \Theta_{\varepsilon
}^{2}\ast f\right)  \}\rangle.
\]
As $f-\gamma_{\varepsilon}^{2}\left(  \Theta_{\varepsilon}^{2}\ast f\right)  $
is a representative of $\sigma(f)-i_{S}(f)$, which is equal to zero in
${\mathcal{G}}(Y)$ since $i_{s}\mid_{{\mathcal{C}}^{\infty}(Y)}=\sigma$, this
representative is in ${\mathcal{I}}(Y)$. Furthermore, $(\gamma_{\varepsilon
}^{1}\Theta_{\varepsilon}^{1}(x-\cdot))_{\varepsilon}$ is in ${\mathcal{E}%
}_{M}(X)$. Thus, $(\{\xi\mapsto\gamma_{\varepsilon}^{1}(\xi)\Theta
_{\varepsilon}^{1}(x-\xi)\}\otimes\{f-\gamma_{\varepsilon}^{2}\left(
\Theta_{\varepsilon}^{2}\ast f\right)  \})_{\varepsilon}$ is in ${\mathcal{I}%
}(X\times Y)$. Let $K_{1}$ be a compact of $X$, then there are
$\Omega_{1}$ relatively compact open subset of $X$ such that
$K_{1}\subset\Omega_{1}$ and $\varepsilon_{1}>0$ such that, for
all $x$ in $K_{1}$, for all $\varepsilon <\varepsilon_{1}$,
$\gamma_{\varepsilon}^{1}\Theta_{\varepsilon}^{1}(x-\cdot)$ is in
${\mathcal{D}}(\Omega_{1})$. Furthermore, by setting
$supp~f=K_{2}$, there are $\Omega_{2}$ relatively compact open
subset of $Y$ such that $K_{2}\subset\Omega_{2}$ and
$\varepsilon_{2}>0$ such that, for all
$\varepsilon<\varepsilon_{2}$, $f-\gamma_{\varepsilon}^{2}\left(
\Theta_{\varepsilon}^{2}\ast f\right)  $ is in
${\mathcal{D}}(\Omega_{2})$. Consequently,
$\{\xi\mapsto\gamma_{\varepsilon}^{1}(\xi)\Theta_{\varepsilon
}^{1}(x-\xi)\}\otimes\{f-\gamma_{\varepsilon}^{2}\left(
\Theta_{\varepsilon }^{2}\ast f\right)  \}$ is in
${\mathcal{D}}(\Omega_{1}\times\Omega_{2})$. By using the local
structure of distributions, one can write $h$ as a derivative of a
continuous function on $\mathbb{R}^{m}\times\mathbb{R}^{n}$, whose
support is contained in an arbitrary neighborhood of
$\Omega_{1}\times \Omega_{2}$ and one shows that
\[
\left(  \langle h,\{\xi\mapsto\gamma_{\varepsilon}^{1}(\xi)\Theta
_{\varepsilon}^{1}(x-\xi)\}\otimes\{f-\gamma_{\varepsilon}^{2}\left(
\Theta_{\varepsilon}^{2}\ast f\right)  \}\rangle\right)  _{\varepsilon}%
\]
is in ${\mathcal{I}}(X)$, which implies the required result.

\begin{proposition}
If $H$ is in ${\mathcal{G}}^{\infty}(X\times Y)$ then $H$ is regular in the
sense given by definition \ref{defreg}.
\end{proposition}

\underline{{\bf{Proof}}}. Let $H$ be in ${\mathcal{G}}^{\infty}(X\times Y)$ and $f$ be in ${\mathcal{G}%
}_{C}(Y)$. There exists $K_{2}$ compact of $Y$ and $(f_{\varepsilon
})_{\varepsilon}$ a representative of $f$ such that $supp~f_{\varepsilon
}\subset K_{2}$. Let us denote by $(H_{\varepsilon})_{\varepsilon}$ a
representative of $H$. A representative of $\widehat{H}(f)$ is
\[
\left(  x,y\right)  \mapsto\int_{K_{2}}H_{\varepsilon}(x,y)f_{\varepsilon
}(y)\,\mathrm{d}y
\]
and, for all $K_{1}$ compact of $X$, $\alpha$ in $\mathbb{N}^{m}$, $x$ in
$K_{1}$,
\[
\left\vert \partial_{x}^{\alpha}\left(  \int_{K_{2}}H_{\varepsilon
}(x,y)f_{\varepsilon}(y)\,\mathrm{d}y\right)  \right\vert \leq Vol(K_{2}%
)\sup_{(x,y)\in K_{1}\times K_{2}}|\partial_{x}^{\alpha}H_{\varepsilon
}(x,y)|\sup_{y\in K_{2}}|f_{\varepsilon}(y)|\leq C(\alpha)\varepsilon^{-q},
\]
as $\varepsilon$ tends to 0, where $q$ does not depend on
$\alpha$. That shows that $\widehat{H}(f)$ is in
${\mathcal{G}}^{\infty}(X)$.\medskip

This result has also been proved in \cite{garetto2}.

\section{Composition of generalized integral operators}

\subsection{Operators with kernel in $\mathcal{G}$$_{ps}\left(  .\right)  $}

\begin{theorem}
\label{ThmCompG}Let $X$, $Y$ and $\Xi$ be three open subsets of $\mathbb{R}%
^{m}$, $\mathbb{R}^{n}$ and $\mathbb{R}^{p}$ respectively and $H_{1}\in$
$\mathcal{G}$$_{ps}(X\times\Xi)$, $H_{2}\in$ $\mathcal{G}$$_{ps}(\Xi\times
Y)$. The operators $\widehat{H}_{1}~:~\mathcal{G}$$(\Xi)\rightarrow
\mathcal{G}$$(X)$ and $\widehat{H}_{2}$~:~$\mathcal{G}$$(Y)\rightarrow
\mathcal{G}$$(\Xi)$ can be composed. Moreover, $\widehat{H}_{1}\circ
\widehat{H}_{2}$ is a generalized integral operator, whose kernel is $L$,
defined by $L(\cdot_{1},\cdot_{2})=\int_{\Xi}H_{1}(\cdot_{1},\xi)H_{2}%
(\xi,\cdot_{2})\,\mathrm{d}\xi$ (with the meaning of notation \ref{NotIntPar})
and $L$ belongs to $\mathcal{G}$$_{ps}(X\times Y)$.
\end{theorem}

\underline{{\bf{Proof}}}. For all $f$ in $\mathcal{G}$$(Y)$,
$\widehat{H}_{2}\left(  f\right)  $ is well
defined in $\mathcal{G}$$(\Xi)$, then we can define $\widehat{H}_{1}%
(\widehat{H}_{2}\left(  f\right)  )$ in $\mathcal{G}$$(X)$ and $\widehat
{H}_{1}\circ\widehat{H}_{2}$ is well defined. We have to show that $L$ is well
defined, properly supported and that $\widehat{H}_{1}\circ\widehat{H}%
_{2}=\widehat{L}$. We set $\Phi(\cdot_{1},\cdot_{2},\cdot_{3})=H_{1}(\cdot
_{1},\cdot_{3})H_{2}(\cdot_{3},\cdot_{2})$ ($\cdot_{3}$ refers to the $\xi$
variable). Choose $O_{1}$ (resp. $O_{2}$) a relatively compact open subset of
$X$ (resp. $Y$).

Since $H_{1}\in$ $\mathcal{G}$$_{ps}(X\times\Xi),$ there exists a
compact subset $K_{3}\subset\Xi$ such that
\[
supp~H_{1}\cap(O_{1}\times\Xi)\subset O_{1}\times K_{3}.
\]
Therefore, $supp~\Phi\cap(O_{1}\times O_{2}\times\Xi)\subset O_{1}\times
O_{2}\times K_{3}$ and $\Phi$ is in $\mathcal{G}$$_{ps}(X\times Y\times\Xi)$.
Proposition \ref{intgps} implies the existence of $L$ in $\mathcal{G}%
$$(X\times Y)$, denoted by $\int H_{1}(\cdot_{1},\xi)H_{2}(\xi,\cdot
_{2})\,\mathrm{d}\xi$.

With the same notations as above, since $H_{2}\in$ $\mathcal{G}$$_{ps}%
(X\times\Xi)$, there exists for $O_{3}=\mathring{K}_{3}$, a
compact subset
$K_{2}\subset Y$ such that%
\[
supp~H_{2}\cap(O_{3}\times Y)\subset O_{3}\times K_{2}.
\]
We have $\Phi_{\left\vert O_{1}\times\left(  Y\backslash K_{2}\right)
\times\Xi\right.  }=0$ since $H_{1\left\vert O_{1}\times\left(  \Xi\backslash
K_{3}\right)  \right.  }=0$ and $H_{2\left\vert K_{3}\times\left(  Y\backslash
K_{2}\right)  \right.  }=0$. Therefore, $L$ is null on $O_{1}\times\left(
Y\backslash K_{2}\right)  $, so $supp~L\cap(O_{1}\times Y)$ is included in
$O_{1}\times K_{2}$ which shows that $L$ is in $\mathcal{G}$$_{ps}(X\times Y)$.

Moreover, for any $f\in\mathcal{G}$$(Y)$, we have (the compact sets on which
the integration are performed are indicated contrary to the notations above),
\begin{align*}
\widehat{L}\left(  f\right)  _{\left\vert O_{1}\right.  }  &  =\int_{K_{2}%
}L\left(  \cdot_{1},y\right)  f\left(  y\right)  \,\mathrm{d}y\\
&  =\int_{K_{2}}\left(  \int_{K_{3}}H_{1}(\cdot_{1},\xi)H_{2}(\xi
,y)\,\mathrm{d}\xi\right)  f\left(  y\right)  \mathrm{d}y\\
&  =\int_{K_{2}\times K_{3}}H_{1}(\cdot_{1},\xi)H_{2}(\xi,y)f\left(  y\right)
\mathrm{d}\xi\mathrm{d}y\ (\text{by\ fubini's\ theorem})\\
&  =\int_{K_{3}}H_{1}(\cdot_{1},\xi)\left(  \int_{K_{2}}H_{2}(\xi,y)f\left(
y\right)  \mathrm{\,\mathrm{d}}\xi\right)  \mathrm{\,\mathrm{d}}y=\widehat
{H}_{1}(\widehat{H}_{2}\left(  f\right)  )_{\left\vert O_{1}\right.  }.
\end{align*}
Using the sheaf structure of $\mathcal{G}$$(X)$, it follows that
$\widehat {L}\left(  f\right)  =\left(
\widehat{H}_{1}\circ\widehat{H}_{2}\right)  (f)$.

\begin{remark}
It is straightforward to verify that the composition of generalized integral
operators is associative.
\end{remark}

\begin{example}
Take $\delta\in{\mathcal{D}}^{\prime}(\mathbb{R})$ the classical delta
function and $\chi$ an integrable function on $\mathbb{R}$. Set $H=\delta
_{x}\otimes\mathbf{1}_{y}$ and $K=\chi\otimes\delta_{y}$. Such $H$ and $K$ are
distributions on $\mathbb{R}^{2}$. The kernel operators associated with $H$
and $K$ are respectively
\[%
\begin{array}
[c]{ccccccccccc}%
\widehat{H} & : & {\mathcal{D}}(\mathbb{R}) & \rightarrow & {\mathcal{D}%
}^{\prime}(\mathbb{R}) & \ \ \ \ \ \mathrm{and}\ \ \ \ \  & \widehat{K} & : &
{\mathcal{D}}(\mathbb{R}) & \rightarrow & {\mathcal{D}}^{\prime}(\mathbb{R})\\
&  & f & \mapsto & \delta_{x}\int f(y)\,\mathrm{d}y &  &  &  & f & \mapsto &
f(0)\chi.
\end{array}
\]
By noticing that $\widehat{H}:L^{1}(\mathbb{R})\rightarrow{\mathcal{D}%
}^{\prime}(\mathbb{R})$ and $\widehat{K}:{\mathcal{D}}(\mathbb{R})\rightarrow
L^{1}(\mathbb{R})$, one can define $\widehat{H}\circ\widehat{K}$ by the
following
\[
\widehat{H}\circ\widehat{K}:%
\begin{array}
[t]{ccc}%
{\mathcal{D}}(\mathbb{R}) & \rightarrow & {\mathcal{D}}^{\prime}(\mathbb{R})\\
f & \mapsto & \delta_{x}f(0)\int\chi(\xi)\,\mathrm{d}\xi,
\end{array}
\]
which admits as kernel
\begin{equation}
\delta_{x}\delta_{y}\int\chi(\xi)\,\mathrm{d}\xi. \label{exnoyau}%
\end{equation}
Conversely, $\widehat{K}\circ\widehat{H}$ cannot be defined classically. We
are going to define it in the context of generalized functions. By using the
notations of paragraph \ref{injection}, one has
\[
i_{S}(H)=Cl\left(  (x,y)\mapsto\left(  \Theta_{\varepsilon}\otimes
\mathbf{1}_{y}\right)  (x,y)\right)  _{\varepsilon}%
\]
and
\[
i_{S}(K)=Cl\left(  (x,y)\mapsto\left(  (\chi\ast\Theta_{\varepsilon}%
)\otimes\Theta_{\varepsilon}\right)  (x,y)\right)  _{\varepsilon},
\]
so $\widehat{i_{S}(K)}\circ\widehat{i_{S}(H)}=\widehat{L}$ is well defined
from ${\mathcal{G}}_{C}(\mathbb{R})$ to ${\mathcal{G}}(\mathbb{R})$, with
\begin{align*}
L  &  =Cl\left(  (x,y)\mapsto\int((\chi\ast\Theta_{\varepsilon})\otimes
\rho_{\varepsilon})(x,\xi)(\Theta_{\varepsilon}\otimes\mathbf{1}_{y}%
)(\xi,y)\,\mathrm{d}\xi\right)  _{\varepsilon}\\
&  =Cl\left(  (x,y)\mapsto(\chi\ast\Theta_{\varepsilon})(x)\int(\Theta
_{\varepsilon}(\xi))^{2}\,\mathrm{d}\xi\right)  _{\varepsilon},
\end{align*}
that is $L=\int i_{S}(\delta)^{2}(\xi)\,\mathrm{d}\xi\cdot i_{S}(\chi
)\otimes\mathbf{1}_{y}=\int\delta^{2}(\xi)\,\mathrm{d}\xi\cdot\chi
\otimes\mathbf{1}_{y}$ with a slight abuse of notation. For the case of
$\widehat{H}\circ\widehat{K}$, one can easily verify that the image by $i_{S}$
of the classical distributional kernel given by (\ref{exnoyau}) is equal to
the kernel obtained by theorem \ref{ThmCompG}.

\end{example}

\begin{corollary}
\label{PropCompItGPS}For $H$ in ${\mathcal{G}}_{ps}(X^{2})$ ($X$ open subset
of $\mathbb{R}^{d}$) and $n\geq2$, $\widehat{H}^{n}=\underset{n\ times}%
{\underbrace{\widehat{H}\circ\cdots\circ\widehat{H}}}:{\mathcal{G}%
}(X)\rightarrow{\mathcal{G}}(X)$ is a well defined generalized integral
operator, whose kernel $L_{n}\in\mathcal{G}$$_{ps}(X^{2})$ is defined by
\[
L_{n}(\cdot_{1},\cdot_{2})=\int H(\cdot_{1},\xi_{1})H(\xi_{1},\xi_{2})\cdots
H(\xi_{n-1},\cdot_{2})\,\mathrm{d}\xi_{1}\mathrm{d}\xi_{2}\cdots\mathrm{d}%
\xi_{n-1},
\]
with the meaning of notation \ref{NotIntPar}.
\end{corollary}

\underline{{\bf{Proof}}}. We prove this proposition by induction.\
Theorem \ref{ThmCompG} gives the result for $n=2$, by considering
$X=\Xi=Y$ and $H_{1}=H_{2}$. Suppose now that $\widehat{H}^{n-1}$
is well defined with its kernel $L_{n-1}$ defined by
\[
L_{n-1}(\cdot_{1},\cdot_{2})=\int_{X^{n-2}}H(\cdot_{1},\xi_{1})H(\xi_{1}%
,\xi_{2})\cdots H(\xi_{n-2},\cdot_{2})\,\mathrm{d}\xi_{1}\mathrm{d}\xi
_{2}\cdots\mathrm{d}\xi_{n-2}%
\]
in ${\mathcal{G}}_{ps}(X^{2})$.

We apply theorem \ref{ThmCompG} with $H_{1}=H$ and $H_{2}=L_{n-1}$.\ It
follows that $\widehat{H}^{n}=\widehat{H}^{n-1}\circ\widehat{H}$ is a well
defined operator which admits as kernel
\begin{align*}
L(\cdot_{1},\cdot_{2})  &  =\int_{X}\left(  \int_{X^{n-2}}H(\cdot_{1},\xi
_{1})H(\xi_{1},\xi_{2})\cdots H(\xi_{n-2},\xi_{n-1})\,\mathrm{d}\xi
_{1}\mathrm{d}\xi_{2}\cdots\mathrm{d}\xi_{n-2}\right)  H(\xi_{n-1},\cdot
_{2})\,\mathrm{d}\xi_{n-1}\\
&  =\int_{X^{n-1}}H(\cdot_{1},\xi_{1})H(\xi_{1},\xi_{2})\cdots H(\xi_{n-2}%
,\xi_{n-1})H(\xi_{n-1},\cdot_{2})\,\mathrm{d}\xi_{1}\mathrm{d}\xi_{2}%
\cdots\mathrm{d}\xi_{n-2}\mathrm{d}\xi_{n-1},
\end{align*}
by applying Fubini's theorem. This is possible since the integrals
are always performed on compact sets, using the local definition
(proposition \ref{intgps}). Thus $L(\cdot_{1},\cdot_{2})$, which
is properly supported, is equal to $L_{n}(\cdot_{1},\cdot_{2})$
which satisfies the required properties.\medskip

For operators with compactly supported kernels, we can give a more
precise result.

\begin{proposition}
\label{corol}With the notations of theorem \ref{ThmCompG}, for $H_{1}$ in
$\mathcal{G}_{C}(X\times\Xi)$ and $H_{2}$ in $\mathcal{G}_{C}(\Xi\times Y)$,
$\widehat{H}_{1}\circ\widehat{H}_{2}:\mathcal{G}(Y)\rightarrow\mathcal{G}%
_{C}(X)$ is a generalized integral operator whose kernel $L$ is an element of
$\mathcal{G}$$_{C}(X\times Y)$.\ Moreover, if $K_{1}$ (resp. $K_{2}~$; $K_{3}%
$) is a compact subset of $X$, (resp. $\Xi~;Y$) such that the
support of $H_{1}$ (resp. $H_{2}$) is contained in the interior of
$K_{1}\times K_{2}$ (resp.
$K_{2}\times K_{3}$) then $L$ can be defined globally by $L(\cdot_{1}%
,\cdot_{2})=\int_{K_{2}}H_{1}(\cdot_{1},\xi)H_{2}(\xi,\cdot_{2})\,\mathrm{d}%
\xi$ and the support of $L$ is contained in $K_{1}\times K_{3}$.
\end{proposition}

\underline{{\bf{Proof}}}. We only have to verify the assertions
related to $L$. Denote by $\left( H_{1,\varepsilon}\right)
_{\varepsilon}$ (resp. $\left(  H_{2,\varepsilon }\right)
$)$_{\varepsilon}$ a representative of $H_{1}$ (resp. $H_{2}$) and
set $O_{1}=X\setminus K_{1}$, $O_{3}=Y\setminus K_{3}$. The net
\[
\left(  \left(  x,y\right)  \mapsto\int_{K_{2}}H_{1,\varepsilon}%
(x,\xi)H_{2,\varepsilon}(\xi,y)\,\mathrm{d}\xi\right)  _{\varepsilon}%
\]
is a representative of $L$, which justifies the global definition of $L$. For
$U\Subset X$ and $V\Subset Y$ such that $U\times V\subset X\times Y\backslash
K_{1}\times K_{3}$, we have either $U\subset O_{1}$ or $V\subset O_{3}$. We
shall suppose, for example, that $U\subset O_{1}$. For $(x,y)\in U\times V$,
we have%
\[
\left.  \left\vert L_{\varepsilon}(x,y)\right\vert =\left\vert \int_{K_{2}%
}H_{1,\varepsilon}(x,\xi)H_{2,\varepsilon}(\xi,y)\,\mathrm{d}\xi\right\vert
\right.  \leq Vol(K_{2})p_{U\times K_{2}}(H_{1,\varepsilon})p_{K_{2}\times
V}(H_{2,\varepsilon}).
\]
Therefore
\begin{equation}
p_{U\times V}(L_{\varepsilon})\leq Vol(K_{2})p_{U\times K_{2}}%
(H_{1,\varepsilon})p_{K_{2}\times V}(H_{2,\varepsilon}). \label{GF04BCD11}%
\end{equation}
As $(H_{1,\varepsilon\left\vert O_{1}\times\Xi\right.  })_{\varepsilon}$ is in
$\mathcal{I}(O_{1}\times\Xi)$ and $U\cap K_{2}\subset O_{1}\times\Xi$, it
follows that $p_{U\times K_{2}}(H_{1,\varepsilon})=O\left(  \varepsilon
^{m}\right)  $ as $\varepsilon\rightarrow0$, for all $m\in\mathbb{N}$. Since
$\left(  H_{2,\varepsilon}\right)  _{\varepsilon}$ is in $\mathcal{E}%
_{M}\left(  \Xi\times Y\right)  $, relation (\ref{GF04BCD11})
implies that $p_{U\times V}(L_{\varepsilon})=O\left(
\varepsilon^{m}\right)  $ as $\varepsilon\rightarrow0$, for all
$m\in\mathbb{N}$. Finally, $(L_{\varepsilon })_{\varepsilon}$
satisfies the null estimate of order $0$ for all compact subsets
included in $X\times Y\backslash K_{1}\times K_{3}$. Therefore,
$L_{\left\vert X\times Y\backslash K_{1}\times K_{3}\right.  }=0$
and the support of $L$ is contained in $K_{1}\times
K_{3}$.\medskip

From corollary \ref{PropCompItGPS} and proposition \ref{corol}, we
immediately deduce the following:

\begin{corollary}
\label{PropCompItG}If $H$ belongs to ${\mathcal{G}}_{C}(X^{2})$,
with $supp~H\subset(\mathring{K}_{2})^{2}$ ($K_{2}\Subset X$),
then the image of $\widehat{H}^{n}$ is included in
${\mathcal{G}}_{C}(X)$ and the support of $L_{n}$ is contained in
$(\mathring{K}_{2})^{2}$. Moreover, $L_{n}\in$
${\mathcal{G}}_{C}(X^{2})$ can be defined globally by
\[
L_{n}(\cdot_{1},\cdot_{2})=\int_{K_{2}}H(\cdot_{1},\xi_{1})H(\xi_{1},\xi
_{2})\cdots H(\xi_{n-1},\cdot_{2})\,\mathrm{d}\xi_{1}\mathrm{d}\xi_{2}%
\cdots\mathrm{d}\xi_{n-1}.
\]

\end{corollary}

\subsection{Operators with kernel in ${\mathcal{G}}_{L^{2}}\left(
\cdot\right)  $}

\begin{proposition}
\label{PropCompGL}For $H_{1}$ in ${\mathcal{G}}_{L^{2}}(X\times\Xi)$ and
$H_{2}$ in ${\mathcal{G}}_{L^{2}}(\Xi\times Y)$, $\widehat{H}_{1}\circ
\widehat{H}_{2}:{\mathcal{G}}_{L^{2}}(Y)\rightarrow{\mathcal{G}}_{L^{2}}(X)$
is a generalized integral operator whose kernel is $L\in\mathcal{G}$$_{L^{2}%
}(X\times Y)$ defined by
\[
L=Cl\left(  \left(  x,y\right)  \mapsto\int_{\Xi}H_{1,\varepsilon}%
(x,\xi)H_{2,\varepsilon}(\xi,y)\,\mathrm{d}\xi\right)  _{\varepsilon},
\]
where $(H_{1,\varepsilon})_{\varepsilon}$ (resp. $(H_{2,\varepsilon
})_{\varepsilon}$) is a representative of $H_{1}$ (resp. $H_{2}$).
\end{proposition}

\underline{{\bf{Proof}}}. With the notations given in the proposition, set%
\[
L_{\varepsilon}(x,y)=\int_{\Xi}H_{1,\varepsilon}(x,\xi)H_{2,\varepsilon}%
(\xi,y)\,\mathrm{d}\xi,\mbox{ for all }(x,y)\mbox{
in }X\times Y,
\]
Then
\begin{align*}
\Vert L_{\varepsilon}\Vert_{2}^{2}  &  =\int\int\left(  \int H_{1,\varepsilon
}(x,\xi)H_{2,\varepsilon}(\xi,y)\,\mathrm{d}\xi\right)  ^{2}\,\mathrm{d}%
x\mathrm{d}y\\
&  \leq\int\int\Vert H_{1,\varepsilon}(x,\cdot)\Vert_{2}^{2}\Vert
H_{2,\varepsilon}(\cdot,y)\Vert_{2}^{2}\,\mathrm{d}x\mathrm{d}y\\
&  \leq\int\Vert H_{1,\varepsilon}(x,\cdot)\Vert_{2}^{2}\,\mathrm{d}x\int\Vert
H_{2,\varepsilon}(\cdot,y)\Vert_{2}^{2}\,\mathrm{d}y=\Vert H_{1,\varepsilon
}\Vert_{2}^{2}\Vert H_{2,\varepsilon}\Vert_{2}^{2}.
\end{align*}
Furthermore, for all $\alpha,\beta\in\mathbb{N}^{d}\setminus\{(0,0)\}$ and
$(x,y)\in X\times Y$, by derivation in the sense of distributions, we have
\[
\partial_{x\,y}^{\left(  \alpha,\beta\right)  }L_{\varepsilon}(x,y)=\int_{\Xi
}\partial_{x}^{\alpha}H_{1,\varepsilon}\left(  x,\xi\right)  \partial
_{y}^{\beta}H_{2,\varepsilon}\left(  \xi,y\right)  \,\mathrm{d}\xi
\ \ \ \ (\mathrm{with}\ \partial_{x\,y}^{\left(  \alpha,\beta\right)
}L_{\varepsilon}=\frac{\partial^{\left\vert \alpha\right\vert +\left\vert
\beta\right\vert }L_{\varepsilon}}{\partial x^{\alpha}\partial y^{\beta}}).
\]
Thus
\[
\left\Vert \partial_{x\,y}^{\left(  \alpha,\beta\right)  }L_{\varepsilon
}\right\Vert _{2}^{2}\leq\left\Vert \partial_{x}^{\alpha}H_{1,\varepsilon
}\right\Vert _{2}^{2}\left\Vert \partial_{y}^{\beta}H_{2,\varepsilon
}\right\Vert _{2}^{2}.
\]
As $(H_{1,\varepsilon})_{\varepsilon}\in{\mathcal{E}}_{L^{2}}(X\times\Xi)$ and
$(H_{2,\varepsilon})_{\varepsilon}\in{\mathcal{E}}_{L^{2}}(\Xi\times Y)$, we
get that $(L_{\varepsilon})_{\varepsilon}$ is in ${\mathcal{E}}_{L^{2}%
}(X\times Y)$. We set $L=Cl(L_{\varepsilon})_{\varepsilon}$ in ${\mathcal{G}%
}_{L^{2}}(X\times Y)$. Moreover, for any $f\in$ ${\mathcal{G}}_{L^{2}}(X)$ and
any of its representative $(f_{\varepsilon})_{\varepsilon}\in\mathcal{E}%
$$_{L^{2}}(\Omega)$, a representative of $\left(  \widehat{H}_{1}\circ
\widehat{H}_{2}\right)  (f)$ is given by the net $\left(  \Psi_{\varepsilon
}\right)  _{\varepsilon}$ with
\[
\Psi_{\varepsilon}(x)=\int_{\Xi}H_{1,\varepsilon}(x,\xi)\left(  \int
_{Y}H_{2,\varepsilon}(\xi,y)f_{\varepsilon}(y)\,\mathrm{d}y\right)
\,\mathrm{d}\xi,\ \ \mathrm{for\ all}\ x\in X.
\]
Then, for all $x\in X,$
\[
\Psi_{\varepsilon}(x)=\int_{\Xi}\int_{Y}H_{1,\varepsilon}(x,\xi
)H_{2,\varepsilon}(\xi,y)f_{\varepsilon}(y)\,\mathrm{d}y\mathrm{d}\xi=\int
_{Y}L_{\varepsilon}(x,y)f_{\varepsilon}(y)\,\mathrm{d}y,
\]
by applying Fubini's theorem. Thus, $\left(  \Psi_{\varepsilon}\right)
_{\varepsilon}$ is a representative of $\widehat{L}(f)$ and $\widehat{H}%
_{1}\circ\widehat{H}_{2}=\widehat{L}$.

\begin{corollary}
\label{PropCompGLIT}For $H$ in ${\mathcal{G}}_{L^{2}}(X^{2})$ ($X$ open subset
of $\mathbb{R}^{d}$), for all $n\geq2$, $\widehat{H}^{n}:{\mathcal{G}}_{L^{2}%
}(X)\rightarrow{\mathcal{G}}_{L^{2}}(X)$ is a generalized integral operator
whose kernel is $L_{n}\in$ ${\mathcal{G}}_{L^{2}}(X^{2})$ defined by
$L_{n}=Cl\left(  L_{n,\varepsilon}\right)  _{\varepsilon}$, with
\[
L_{n,\varepsilon}:\left(  x,y\right)  \mapsto\int_{X^{n-1}}H_{\varepsilon
}(x,\xi_{1})H_{\varepsilon}(\xi_{1},\xi_{2})\cdots H_{\varepsilon}(\xi
_{n-1},y)\,\mathrm{d}\xi_{1}\mathrm{d}\xi_{2}\cdots\mathrm{d}\xi_{n-1},
\]
where $(H_{\varepsilon})_{\varepsilon}$ is a representative of $H$.
\end{corollary}

\underline{{\bf{Proof}}}. The proof of this result is an
adaptation of the one of corollary \ref{PropCompItGPS} to the
$L^2$-case.

\section{Application: Exponential of generalized integral
operators}\label{exponentiel}

In this section, we define the exponential of generalized integral
operators in two particular cases and study some of their
properties. We first define some convenient spaces. Let us set,
for $\Omega$ open subset of $\mathbb{R}^{d}$ ($d\in\mathbb{N}$),
\[
{\mathcal{H}}_{ln}(\Omega)=\left\{  (u_{\varepsilon})_{\varepsilon}%
\in{\mathcal{C}}^{\infty}(\Omega)^{(0,1]}~/~\forall K\Subset\Omega~,~\forall
l\in\mathbb{N}~,~p_{K,l}(u_{\varepsilon})=O(|\ln\varepsilon
|)\mbox{ as }\varepsilon\rightarrow0\right\}
\]
and
\[
{\mathcal{H}}_{L^{2}\,ln}(\Omega)=\left\{  (u_{\varepsilon})_{\varepsilon}\in
H^{\infty}(\Omega)^{(0,1]}~/~\forall m\geq0~,~\Vert u_{\varepsilon}\Vert
_{m}=O(|\ln\varepsilon|)\mbox{ as
}\varepsilon\rightarrow0\right\}  .
\]
Then ${\mathcal{H}}_{ln}(\Omega)$ (resp. ${\mathcal{H}}_{L^{2}\,ln}(\Omega)$)
is a linear subspace of ${\mathcal{E}}_{M}(\Omega)$ (resp. ${\mathcal{E}%
}_{L^{2}}(\Omega)$), but is not a subalgebra. Then we set
\[
{\mathcal{G}}_{ln}(\Omega)={\mathcal{H}}_{ln}(\Omega)/{\mathcal{I}}%
(\Omega)\,,\ \ {\mathcal{G}}_{C\,ln}(\Omega)={\mathcal{G}}_{ln}(\Omega
)\cap{\mathcal{G}}_{C}(\Omega)\ \ \ \mathrm{and}\ \ \ {\mathcal{G}}%
_{L^{2}\,ln}(\Omega)={\mathcal{H}}_{L^{2}\,ln}(\Omega))/{\mathcal{I}}_{L^{2}%
}(\Omega).
\]
The space ${\mathcal{G}}_{ln}(\Omega)$ (resp. ${\mathcal{G}}_{L^{2}%
\,ln}(\Omega)$) is a subvector space of ${\mathcal{G}}(\Omega)$ (resp.
${\mathcal{G}}_{L^{2}}(\Omega)$).

\subsection{Exponential of generalized integral operators whose kernel is in
${\mathcal{G}}_{C\,ln}(X^{2})$\label{SbSExpComp}}

\begin{theoremdef}
\label{defexp1}Let $H$ be in $\mathcal{G}_{C\,ln}(X^{2})$ ($X$ an
open subset of $\mathbb{R}^{d}$). For $n\geq1$, denote by $L_{n}$
the kernel of $\widehat
{H}^{n}:\mathcal{G}(X)\rightarrow\mathcal{G}_{C}(X)$ defined as in
corollary \ref{PropCompItG} (with $L_{1}=H$) and by $\left(
L_{n,\varepsilon}\right) _{\varepsilon}$ a representative of
$L_{n}$. For all $\varepsilon\in(0,1]$, the series
$\sum_{n\geq1}\frac{L_{n,\varepsilon}}{n!}$ converges for the
usual topology of ${\mathcal{C}}^{\infty}(X^{2})$.\ Its sum,
denoted by $S_{\varepsilon}$, defines an element $\left(
S_{\varepsilon}\right) _{\varepsilon}$ of
${\mathcal{E}}_{M}(X^{2})$. Furthermore, $S=Cl\left(
S_{\varepsilon}\right)  _{\varepsilon}$ defines a compactly
supported element of $\mathcal{G}(X^{2})$ only depending on
$H$.\newline The well defined operator
$e^{\widehat{H}}=\widehat{S}+Id$ (where $Id$ is the operator
identity) will be called the \emph{exponential} of $\widehat{H}$.
\end{theoremdef}

\underline{{\bf{Proof}}}. We divide it in three parts. The first
part contains the estimates of
$\sum_{n\geq1}\frac{L_{n,\varepsilon}}{n!}$ for a particular
representative of $L_{n}$, constructed from a representative
$\left(  H_{\varepsilon}\right) _{\varepsilon}$. The second part
deals with the independence of $Cl\left( S_{\varepsilon}\right)
_{\varepsilon}$ of the chosen representative of $H$. The third
part shows that $S$ is compactly supported.\smallskip

$\bullet$ Let $H$ be in $\mathcal{G}_{C\,ln}(X^{2})$ and $\left(
H_{\varepsilon}\right) _{\varepsilon}$ one of its representative.
According
to corollary \ref{PropCompItG}, we have $\widehat{H}^{n}=\widehat{L_{n}%
}:\mathcal{G}(X)\rightarrow\mathcal{G}_{C}(X)$ and $L_{n}\in$ $\mathcal{G}%
_{C}(X^{2})$ admits as representative $\left(  L_{n,\varepsilon}\right)
_{\varepsilon}$ with
\[
L_{n,\varepsilon}(\cdot_{1},\cdot_{2})=\int_{K^{n-1}}H_{\varepsilon}(\cdot
_{1},\xi_{1})H_{\varepsilon}(\xi_{1},\xi_{2})\cdots H_{\varepsilon}(\xi
_{n-1},\cdot_{2})\,\mathrm{d}\xi_{1}\mathrm{d}\xi_{2}\cdots\mathrm{d}\xi
_{n-1},
\]
where $K$ is a compact subset of $X$ such that
$supp$~$H\subset\mathring{K}^{2}$.

For all compact subset of $X^{2}$ of the form $K_{1}\times K_{2}$, $\left(
\alpha,\beta\right)  \in\mathbb{N}^{d}\times\mathbb{N}^{d}$ and $(x,y)\in
K_{1}\times K_{2}$, one has%
\[
\left\vert \partial_{x\,y}^{\left(  \alpha,\beta\right)  }L_{2,\varepsilon
}(x,y)\right\vert =\left\vert \int_{K}\partial_{x}^{\alpha}H_{\varepsilon
}(x,\xi)\partial_{y}^{\beta}H_{\varepsilon}(\xi,y)\,\mathrm{d}\xi\right\vert
\leq\int_{K}p_{K_{1}\times K,|\alpha|}(H_{\varepsilon})p_{K\times K_{2}%
,|\beta|}(H_{\varepsilon})\,\mathrm{d}\xi.
\]
It follows that
\[
p_{K_{1}\times K_{2},|(\alpha,\beta)|}(L_{2,\varepsilon})\leq Vol(K)p_{V^{2}%
,|(\alpha,\beta)|}^{2}(H_{\varepsilon}),
\]
where $V$ is a compact subset of $X$ containing $K$,$K_{1}$ and
$K_{2}$. By induction, we show that, for all $n\geq2$,
\[
p_{K_{1}\times K_{2},|(\alpha,\beta)|}(L_{n,\varepsilon})\leq Vol(K)^{n-1}%
p_{V^{2},|(\alpha,\beta)|}^{n}(H_{\varepsilon}).
\]
This last inequality implies that the series $\sum_{n\geq1}\frac
{L_{n,\varepsilon}}{n!}$ converges for the usual topology of
${\mathcal{C}}^{\infty}\left(X^{2}\right)$. Set, for
$\varepsilon\in(0,1]$,%
\[
S_{\varepsilon}=\sum_{n=1}^{+\infty}\frac{L_{n,\varepsilon}}{n!}.
\]
As $L_{n}$ is in $\mathcal{G}_{C}(X^{2})$ and since the convergence is uniform
on compact sets, $S_{\varepsilon}$ belongs to $\mathcal{C}^{\infty}(X^{2})$
for all $\varepsilon\in(0,1]$. Furthermore, for all compact subset of $X^{2}$
of the form $K_{1}\times K_{2}$ and $\left(  \alpha,\beta\right)
\in\mathbb{N}^{d}\times\mathbb{N}^{d}$, one has
\begin{multline*}
p_{K_{1}\times K_{2},|(\alpha,\beta)|}(S_{\varepsilon})\leq\sum_{n=1}%
^{+\infty}\frac{1}{n!}p_{K_{1}\times K_{2},|(\alpha,\beta)|}(L_{n,\varepsilon
})\\
\leq\sum_{n=1}^{+\infty}\frac{1}{n!}Vol(K)^{n-1}p_{V^{2},|(\alpha,\beta)|}%
^{n}(H_{\varepsilon})\leq\frac{1}{Vol(K)}e^{Vol(K)p_{V^{2},|(\alpha,\beta
)|}(H_{\varepsilon})}.
\end{multline*}
Since $H$ is in $\mathcal{G}_{C\,ln}(X^{2})$, $p_{V^{2},|(\alpha,\beta
)|}(H_{\varepsilon})=O(|\ln\varepsilon|)$ as $\varepsilon\rightarrow0$, that
is there exists $k\in\mathbb{N}$ such that $p_{V^{2},|(\alpha,\beta
)|}(H_{\varepsilon})\leq k\ln\left(  1/\varepsilon\right)  $, so
\[
p_{K_{1}\times K_{2},|(\alpha,\beta)|}(S_{\varepsilon})\leq C_{K}%
\varepsilon^{-kVol(K)},
\]
where $C_{K}$ is a constant depending only on $K$. Consequently,
$(S_{\varepsilon})_{\varepsilon}$ is in $\mathcal{E}_{M}(X^{2})$ and we denote
by $S$ its class in $\mathcal{G}(X^{2})$.\medskip

$\bullet$ Let us show now that $S$ does not depend on the choice
of the representative of $H$. Let
$(H_{\varepsilon}^{1})_{\varepsilon}$ and
$(H_{\varepsilon}^{2})_{\varepsilon}$ be two representatives of
$H$ in $\mathcal{G}_{C\,ln}(X^{2})$. From
$(H_{\varepsilon}^{1})_{\varepsilon}$, we
define $(L_{n,\varepsilon}^{1})_{\varepsilon}$ and $(S_{\varepsilon}%
^{1})_{\varepsilon}$, and from $(H_{\varepsilon}^{2})_{\varepsilon}$,
$(L_{n,\varepsilon}^{2})_{\varepsilon}$ and $(S_{\varepsilon}^{2}%
)_{\varepsilon}$. Let $K$ be a compact subset of $X$ such that the
support of $H$ is contained in the interior of $K^{2}$. For all
$n\geq2$, $K_{1}\times K_{2}$ compact subset of $X^{2}$ and
$(x,y)\in K_{1}\times K_{2}$, one has
\begin{multline*}
\left(  L_{n,\varepsilon}^{1}-L_{n,\varepsilon}^{2}\right)  (x,y)=\int
_{K^{n-1}}H_{\varepsilon}^{1}(x,\xi_{1})H_{\varepsilon}^{1}(\xi_{1},\xi
_{2})\cdots H_{\varepsilon}^{1}(\xi_{n-1},y)\,\mathrm{d}\xi_{1}\mathrm{d}%
\xi_{2}\cdots\mathrm{d}\xi_{n-1}\\
-\int_{K^{n-1}}H_{\varepsilon}^{2}(x,\xi_{1})H_{\varepsilon}^{2}(\xi_{1}%
,\xi_{2})\cdots H_{\varepsilon}^{2}(\xi_{n-1},y)\,\mathrm{d}\xi_{1}%
\mathrm{d}\xi_{2}\cdots\mathrm{d}\xi_{n-1}.
\end{multline*}
Thus, $L_{n,\varepsilon}^{1}-L_{n,\varepsilon}^{2}$ can be written
as a sum of $n$ integrals. The term under each integral sign is
itself formed by the product of $n$ functions, one of them being
equal to $H_{\varepsilon}^{1}-H_{\varepsilon}^{2}$ and the $(n-1)$
others being equal to $H_{\varepsilon}^{1}$ or
$H_{\varepsilon}^{2}$. Consequently,
$$p_{K_{1}\times K_{2}}\left(  L_{n,\varepsilon}^{1}-L_{n,\varepsilon}%
^{2}\right)  \leq n Vol(K)^{n-1} p_{V^{2}}\left(  H_{\varepsilon
}^{1}-H_{\varepsilon}^{2}\right) \left(\max_{1\le i\le 2}
(p_{V^{2}}(H_{\varepsilon} ^{i}))\right)^{n-1}.$$
Since $(H_{\varepsilon}^{1})_{\varepsilon}$ and $(H_{\varepsilon}%
^{2})_{\varepsilon}$ are in $\mathcal{G}_{C\,ln}(X^{2})$, there exists
$k\in\mathbb{N}$ such that
$$\max_{1\le i\le 2}
(p_{V^{2}}(H_{\varepsilon} ^{i}))\leq k|\ln\varepsilon|.$$ As
$(H_{\varepsilon}^{1}-H_{\varepsilon}^{2})_{\varepsilon}$ is in
$\mathcal{I}(X^{2})$, for all $m\in\mathbb{N}$, there exists $C>0$
such that $p_{V^{2}}\left(
H_{\varepsilon}^{1}-H_{\varepsilon}^{2}\right)  \leq
C\varepsilon^{m}$ for $\varepsilon$ small enough. Then
\[
p_{K_{1}\times K_{2}}\left(  L_{n,\varepsilon}^{1}-L_{n,\varepsilon}%
^{2}\right)  \leq C\varepsilon^{m}nVol(K)^{n-1}k^{n-1}|\ln\varepsilon
|^{n-1},\ \text{for }\varepsilon\text{ small enough.}%
\]
Thus $$p_{K_{1}\times K_{2}}\left(
S_{\varepsilon}^{1}-S_{\varepsilon}^{2}\right) \leq
C\varepsilon^{m}\sum_{n=1}^{+\infty}\frac{1}{\left(  n-1\right)
!}Vol(K)^{n-1}k^{n-1}|\ln\varepsilon|^{n-1}\leq C\varepsilon^{m}%
\varepsilon^{-kVol(K)},$$
for $\varepsilon$ small enough.
Therefore, $p_{K_{1}\times K_{2}}\left(
S_{\varepsilon}^{1}-S_{\varepsilon}^{2}\right)  =O\left(
\varepsilon
^{m^{\prime}}\right)  $ as $\varepsilon\rightarrow0$, for all $m^{\prime}%
\in\mathbb{N}$, that is $(S_{\varepsilon}^{1}-S_{\varepsilon}^{2}%
)_{\varepsilon}$ belongs to $\mathcal{I}(X^{2})$. Consequently, $S$ does not
depend on the choice of the representative of $H$ in $\mathcal{G}%
_{C\,ln}(X^{2})$.\medskip

$\bullet$ It remains to prove that $S$ is in $\mathcal{G}%
_{C}(X^{2})$, in order to define $\widehat{S}$. Set
$O=X^{2}\setminus K^{2}$. For all $n\geq2$, $K_{1}\times K_{2}$
compact subset of $O$, one has
\[
p_{K_{1}\times K_{2}}(L_{n,\varepsilon})\leq Vol(K)^{n-1}p_{K_{1}\times
K}(H_{\varepsilon})p_{K^{2}}^{n-2}(H_{\varepsilon})p_{K\times K_{2}%
}(H_{\varepsilon}).
\]
Therefore%
\[
p_{K_{1}\times K_{2}}(S_{\varepsilon})\leq p_{K_{1}\times K_{2}}%
(H_{\varepsilon})+\sum_{n=2}^{+\infty}\frac{1}{n!}Vol(K)^{n-1}p_{K_{1}\times
K}(H_{\varepsilon})p_{K^{2}}^{n-2}(H_{\varepsilon})p_{K\times K_{2}%
}(H_{\varepsilon}).
\]
We have either $K_{1}\cap K=\varnothing$ or $K_{2}\cap K=\varnothing$ since
$K_{1}\times K_{2}\subset O$. Suppose, for example, that $K_{1}\cap
K=\varnothing$. Then $K_{1}\times K\subset O$. As $H$ is in $\mathcal{I}(O)$
there exists, for all $m\in\mathbb{N}$, a constant $C>0$ such that
\[
p_{K_{1}\times K_{2}}(H_{\varepsilon})\leq C\varepsilon^{m}\,,\ \ \ ~p_{K_{1}%
\times K}(H_{\varepsilon})\leq C\varepsilon^{m}.
\]
As $H$ is in ${\mathcal{H}}_{ln}(X^{2})$, there exists $k>0$ such that%
\[
p_{K\times K_{2}}(H_{\varepsilon})\leq k|\ln\varepsilon|\,,\ \ \ p_{K^{2}%
}(H_{\varepsilon})\leq k|\ln\varepsilon|.
\]
Thus%
\begin{align*}
p_{K_{1}\times K_{2}}(S_{\varepsilon})  &  \leq C\varepsilon^{m}%
+C\varepsilon^{m}\sum_{n=2}^{+\infty}\frac{1}{\left(  n-1\right)
!}Vol(K)^{n-1}k^{n-1}|\ln\varepsilon|^{n-1}\\
&  \leq
C\varepsilon^{m}+C\varepsilon^{m}e^{Vol(K)k|\ln\varepsilon|}
=O\left(  \varepsilon^{m-kVol(K)}\right)  \text{ as }\varepsilon
\rightarrow0,
\end{align*}
which implies that $(S_{\varepsilon})_{\varepsilon}$ is in
$\mathcal{I}(O)$. Consequently, $S$ is null on $O$ and has a
compact support included in $K^{2}$.

\subsection{Exponential of generalized integral operators whose kernel is in
${\mathcal{G}}_{L^{2}\,ln}(X^{2})$\label{SbSExpGL2}}

\begin{theoremdef}
\label{defexp2}Let $H$ be in ${\mathcal{G}}_{L^{2}\,ln}(X^{2})$ ($X$ open
subset of $\mathbb{R}^{d}$). Denote by $L_{n}$ the kernel of $\widehat{H}%
^{n}:{\mathcal{G}}_{L^{2}}(X)\rightarrow{\mathcal{G}}_{L^{2}}(X)$
defined as in corollary \ref{PropCompGLIT} (with $L_{1}=H$) and by
$\left( L_{n,\varepsilon}\right)  _{\varepsilon}$ a representative
of $L_{n}$. For
all $\varepsilon\in(0,1]$, the series $\sum_{n\geq1}\frac{L_{n,\varepsilon}%
}{n!}$, where $\left(  L_{n,\varepsilon}\right)  _{\varepsilon}\in
{\mathcal{H}}_{L^{2}\,ln}(X^{2})$ is a representative of $L_{n}$,
converges in $L^{2}$-norm. Its sum, denoted by $S_{\varepsilon}$,
defines an
element $\left(  S_{\varepsilon}\right)  _{\varepsilon}$ of ${\mathcal{E}%
}_{L^{2}}(X^{2})$.\newline By setting $S=Cl\left(  S_{\varepsilon}\right)
_{\varepsilon}$ in ${\mathcal{G}}_{L^{2}}(X^{2})$, we define the
\emph{exponential }of $\widehat{H}$ as $e^{\widehat{H}}=\widehat{S}+Id$ where
$Id$ is the operator identity.
\end{theoremdef}

\underline{{\bf{Proof}}}. We divide it in two parts. The first
part contains the estimates of
$\sum_{n\geq1}\frac{L_{n,\varepsilon}}{n!}$ for a particular
representative of $L_{n}$, constructed from a representative
$\left(  H_{\varepsilon}\right) _{\varepsilon}$ of $H,$ and shows
the existence of $\left(  S_{\varepsilon }\right)
_{\varepsilon}$. The second part deals with the independence of
$Cl\left(  S_{\varepsilon}\right)  _{\varepsilon}$ with respect to
the chosen representative of $H$.\smallskip

$\bullet$ Let $H$ be in ${\mathcal{G}}_{L^{2}\,ln}(X^{2})$.
Applying corollary \ref{PropCompGLIT}, we have
$\widehat{H}^{n}=\widehat
{L_{n}}:{\mathcal{G}}_{L^{2}}(X)\rightarrow{\mathcal{G}}_{L^{2}}(X)$,
for all $n\geq2$, with $L_{n}\in$ ${\mathcal{G}}_{L^{2}}(X^{2})$
defined by $L_{n}=Cl\left(  L_{n,\varepsilon}\right)  $ where, for
all $\varepsilon \in(0,1]$ and $(x,y)\in X^{2}$,
\[
L_{n,\varepsilon}(x,y)=\int H_{\varepsilon}(x,\xi_{1})H_{\varepsilon}(\xi
_{1},\xi_{2})\cdots H_{\varepsilon}(\xi_{n-1},y)\,\mathrm{d}\xi_{1}%
\mathrm{d}\xi_{2}\cdots\mathrm{d}\xi_{n-1}%
\]
and $(H_{\varepsilon})_{\varepsilon}$ denote a representative of
$H$.\\
We are going to show by induction that $\left\Vert
L_{n,\varepsilon}\right\Vert _{2}^{n} \leq\left\Vert
H_{\varepsilon}\right\Vert _{2}^{n}$, for all $n$ greater than
$2$. First, for all $\varepsilon\in(0,1]$ and $(x,y)\in X^{2}$,
one has
\[
\left\vert L_{2,\varepsilon}(x,y)\right\vert =\left\vert \int H_{\varepsilon
}(x,\xi)H_{\varepsilon}(\xi,y)\,\mathrm{d}\xi\right\vert \leq\left\Vert
H_{\varepsilon}(x,\cdot)\right\Vert _{2}\left\Vert H_{\varepsilon}%
(\cdot,y)\right\Vert _{2}.
\]
Thus $\left\Vert L_{2,\varepsilon}\right\Vert _{2}\leq\left\Vert
H_{\varepsilon }\right\Vert _{2}^{2}$ and the first step is done.
Assume that, for all $(x,y)\in X^{2}$,
$$\left\vert L_{n-1,\varepsilon}(x,y)\right\vert \leq\left\Vert H_{\varepsilon
}(x,\cdot)\right\Vert _{2}\Vert
H_{\varepsilon}\Vert_{2}^{n-3}\left\Vert
H_{\varepsilon}(\cdot,y)\right\Vert _{2} \mbox{ and } \left\Vert
L_{n-1,\varepsilon}\right\Vert _{2}\leq\left\Vert H_{\varepsilon
}\right\Vert _{2}^{n-1}.$$ Then, for all $(x,y)\in X^{2}$,
\[
\left\vert L_{n,\varepsilon}(x,y)\right\vert =\left\vert \int
L_{n-1,\varepsilon}(x,\xi)H_{\varepsilon}(\xi,y)\,\mathrm{d}\xi\right\vert
\leq\left\Vert L_{n-1,\varepsilon}(x,\cdot)\right\Vert _{2}\left\Vert
H_{\varepsilon}(\cdot,y)\right\Vert _{2}.
\]
As
$$\left\Vert L_{n-1,\varepsilon}(x,\cdot)\right\Vert _{2}^{2}
\leq\int\left\Vert H_{\varepsilon}(x,\cdot)\right\Vert
_{2}^{2}\Vert H_{\varepsilon}\Vert_{2}^{2n-6}\left\Vert
H_{\varepsilon}(\cdot,y)\right\Vert _{2}^{2}\,\mathrm{d}y
\leq\left\Vert H_{\varepsilon}(x,\cdot)\right\Vert _{2}^{2}\Vert
H_{\varepsilon}\Vert_{2}^{2n-4},$$ we get $\left\Vert
L_{n-1,\varepsilon}(x,\cdot)\right\Vert _{2}\leq\left\Vert
H_{\varepsilon}(x,\cdot)\right\Vert _{2}\Vert
H_{\varepsilon}\Vert_{2}^{n-2}$ and $\left\vert
L_{n,\varepsilon}(x,y)\right\vert \leq\left\Vert
H_{\varepsilon}(x,\cdot)\right\Vert _{2}\Vert H_{\varepsilon}\Vert_{2}%
^{n-2}\left\Vert H_{\varepsilon}(\cdot,y)\right\Vert _{2}.$ Consequently,
$$\left\Vert L_{n,\varepsilon}\right\Vert _{2}^{2}
  \leq\int\int\left\Vert H_{\varepsilon}(x,\cdot)\right\Vert _{2}^{2}\Vert
  H_{\varepsilon}\Vert_{2}^{2n-4}\left\Vert H_{\varepsilon}(\cdot,y)\right\Vert
_{2}^{2}\,\mathrm{d}x\mathrm{d}y =\Vert
H_{\varepsilon}\Vert_2^{2n},$$ that is
\begin{equation}\left\Vert
L_{n,\varepsilon}\right\Vert _{2}\leq\left\Vert H_{\varepsilon
}\right\Vert _{2}^{n}. \label{eq1}%
\end{equation}
With a similar method, we can prove that, for all $n\geq2$, $\left(
\alpha,\beta\right)  \in\mathbb{N}^{d}\times\mathbb{N}^{d}\setminus\{(0,0)\}$
and $(x,y)\in X^{2}$,
\[
\left\vert \partial_{x\,y}^{\left(  \alpha,\beta\right)  }L_{n,\varepsilon
}(x,y)\right\vert \leq\left\Vert \partial_{x}^{\alpha}H_{\varepsilon}%
(x,\cdot)\right\Vert _{2}\Vert H_{\varepsilon}\Vert_{2}^{n-2}\left\Vert
\partial_{y}^{\beta}H_{\varepsilon}(\cdot,y)\right\Vert _{2},
\]
and
\begin{equation}
\left\Vert \partial_{x\,y}^{\left(  \alpha,\beta\right)  }L_{n,\varepsilon
}\right\Vert _{2}\leq\left\Vert \partial_{x}^{\alpha}H_{\varepsilon
}\right\Vert _{2}\Vert H_{\varepsilon}\Vert_{2}^{n-2}\left\Vert \partial
_{y}^{\beta}H_{\varepsilon}\right\Vert _{2}. \label{eq2}%
\end{equation}
From the inequalities (\ref{eq1}) and (\ref{eq2}), we deduce that the series
$\sum_{n\geq1}\frac{L_{n,\varepsilon}}{n!}$ and $\sum_{n\geq1}\frac{1}%
{n!}\partial_{x\,y}^{\left(  \alpha,\beta\right)
}L_{n,\varepsilon}$ converge, in $L^{2}$-norm. We set
\[
S_{\varepsilon}(x,y)=\sum_{n=1}^{+\infty}\frac{L_{n,\varepsilon}}%
{n!}(x,y),\mbox{ for all }(x,y)\in X^{2},
\]
and
\[
D_{\varepsilon}^{\alpha,\beta}(x,y)=\sum_{n=1}^{+\infty}\frac{1}{n!}%
\partial_{x\,y}^{\left(  \alpha,\beta\right)  }L_{n,\varepsilon}%
(x,y),\mbox{ for all }\left(  \alpha,\beta\right)  \in\left(  \mathbb{N}%
^{d}\times\mathbb{N}^{d}\right)  \setminus\{(0,0)\},(x,y)\in X^{2}.
\]
We turn to the study of these series. One has
\[
\left\Vert S_{\varepsilon}\right\Vert _{2}\leq\sum_{n=1}^{+\infty}\frac{1}%
{n!}\left\Vert L_{n,\varepsilon}\right\Vert _{2}\leq\sum_{n=1}^{+\infty}%
\frac{1}{n!}\left\Vert H_{\varepsilon}\right\Vert _{2}^{n}\leq e^{\left\Vert
H_{\varepsilon}\right\Vert _{2}}.
\]
As $H$ is in ${\mathcal{G}}_{L^{2}\,ln}(X^{2})$, $H_{\varepsilon}$ is in
$L^{2}(X^{2})$ and $S_{\varepsilon}$ also, for all $\varepsilon\in(0,1]$.
Furthermore, $\left\Vert H_{\varepsilon}\right\Vert _{2}=O(|\ln\varepsilon|)$
as $\varepsilon$ tends to $0$, that is there exists $k\in\mathbb{N}$ such that
$\left\Vert H_{\varepsilon}\right\Vert _{2}\leq k\ln\left\vert \varepsilon
\right\vert $.\ It follows that
\[
\left\Vert S_{\varepsilon}\right\Vert _{2}\leq\varepsilon^{-k}.
\]
Furthermore, a straightforward exercise on distributions theory
shows that, for all
$\alpha,\beta\in\mathbb{N}^{d}\setminus\{(0,0)\}$ and
$\varepsilon\in(0,1]$, $
\partial_{x\,y}^{\left(  \alpha,\beta\right)  }S_{\varepsilon}=D_{\varepsilon
}^{\alpha,\beta}\mbox{ in }{\mathcal{D}}^{\prime}(X^{2})$.\\
Consequently,
\begin{align*}
\left.  \left\Vert \partial_{x\,y}^{\left(  \alpha,\beta\right)
}S_{\varepsilon}\right\Vert _{2}=\left\Vert D_{\varepsilon}^{\alpha,\beta
}\right\Vert _{2}\right.   &  \leq\sum_{n=1}^{+\infty}\frac{1}{n!}\left\Vert
\partial_{x\,y}^{\left(  \alpha,\beta\right)  }L_{n,\varepsilon}\right\Vert
_{2}\\
&  \leq\left\Vert \partial_{x\,y}^{\left(  \alpha,\beta\right)  }%
H_{\varepsilon}\right\Vert _{2}+\sum_{n=2}^{+\infty}\frac{1}{n!}\left\Vert
\partial_{x}^{\alpha}H_{\varepsilon}\right\Vert _{2}\Vert H_{\varepsilon}%
\Vert_{2}^{n-2}\left\Vert \partial_{y}^{\beta}H_{\varepsilon}\right\Vert
_{2}\\
&  \leq\left\Vert \partial_{x\,y}^{\left(  \alpha,\beta\right)  }%
H_{\varepsilon}\right\Vert _{2}+\left\Vert \partial_{x}^{\alpha}%
H_{\varepsilon}\right\Vert _{2}\left\Vert \partial_{y}^{\beta}H_{\varepsilon
}\right\Vert _{2}e^{\left\Vert H_{\varepsilon}\right\Vert _{2}}.
\end{align*}
As $H$ is in ${\mathcal{G}}_{L^{2}\,ln}(X^{2})$, there exists $k\in\mathbb{N}$
such that $\left\Vert H_{\varepsilon}\right\Vert _{2}$, $\left\Vert
\partial_{x}^{\alpha}H_{\varepsilon}\right\Vert _{2}$, $\left\Vert
\partial_{y}^{\beta}H_{\varepsilon}\right\Vert _{2}$ and $\left\Vert
\partial_{x\,y}^{\left(  \alpha,\beta\right)  }H_{\varepsilon}\right\Vert
_{2}$ are less than $\ln\left(  \varepsilon^{-k}\right)  $, for $\varepsilon$
small enough. Hence
\[
\left\Vert \partial_{x\,y}^{\left(  \alpha,\beta\right)  }S_{\varepsilon
}\right\Vert _{2}=O\left(  \varepsilon^{-k-1}\right)  \text{ as }%
\varepsilon\rightarrow0.
\]
Finally, $(S_{\varepsilon})_{\varepsilon}$ is in ${\mathcal{E}}_{L^{2}}%
(X^{2})$ and we can denote by $S$ its class in ${\mathcal{G}}_{L^{2}}(X^{2}%
)$.\medskip

$\bullet$ We show now that $S$ does not depend on the choice of\
the representative of $H$. Let
$(H_{\varepsilon}^{1})_{\varepsilon}$ and
$(H_{\varepsilon}^{2})_{\varepsilon}$ be two representatives of
$H$ in
${\mathcal{G}}_{L^{2}ln}(X^{2})$. As previously, from $(H_{\varepsilon}%
^{1})_{\varepsilon}$, we define $(L_{n,\varepsilon}^{1})_{\varepsilon}$ and
$(S_{\varepsilon}^{1})_{\varepsilon}$, and from $(H_{\varepsilon}%
^{2})_{\varepsilon}$, we define
$(L_{n,\varepsilon}^{2})_{\varepsilon}$ and
$(S_{\varepsilon}^{2})_{\varepsilon}$. For $n$ greater than 2, we
write $L_{n,\varepsilon}^{1}-L_{n,\varepsilon}^{2}$, as in the
proof of theorem-definition \ref{defexp1}, as a sum of $n$
integrals. After derivating this last expression, we obtain, for
$\left( \alpha ,\beta\right)
\in\mathbb{N}^{d}\times\mathbb{N}^{d}$,

\begin{multline*}
\left\Vert \partial_{x\,y}^{\left(  \alpha,\beta\right)  }\left(
L_{n,\varepsilon}^{1}-L_{n,\varepsilon}^{2}\right)  \right\Vert _{2}%
\leq\left\Vert \partial_{x}^{\alpha}(H_{\varepsilon}^{1}-H_{\varepsilon}%
^{2})\right\Vert _{2}\left\Vert H_{\varepsilon}^{1}\right\Vert _{2}%
^{n-2}\left\Vert \partial_{y}^{\beta}H_{\varepsilon}^{1}(\cdot,y)\right\Vert
_{2}\\
+\left\Vert \partial_{x}^{\alpha}H_{\varepsilon}^{2}\right\Vert _{2}\left\Vert
H_{\varepsilon}^{1}-H_{\varepsilon}^{2}\right\Vert _{2}\left\Vert
H_{\varepsilon}^{1}\right\Vert _{2}^{n-3}\left\Vert \partial_{y}^{\beta
}H_{\varepsilon}^{1}\right\Vert _{2}\\
+\cdots\\
+\left\Vert \partial_{x}^{\alpha}H_{\varepsilon}^{2}\right\Vert _{2}\left\Vert
H_{\varepsilon}^{2}\right\Vert _{2}^{n-2}\left\Vert \partial_{y}^{\beta
}(H_{\varepsilon}^{1}-H_{\varepsilon}^{2})\right\Vert _{2}.
\end{multline*}
Since $(H_{\varepsilon}^{1})_{\varepsilon}$ and $(H_{\varepsilon}%
^{2})_{\varepsilon}$ are in ${\mathcal{G}}_{L^{2}\,ln}(X^{2})$, there exists
$k\in\mathbb{N}$ such that for $\gamma=\left(  0,0\right)  $, $\gamma=\left(
\alpha,0\right)  $ and $\gamma=\left(  0,\beta\right)  $,%
\[
\left\Vert \partial_{x\,y}^{\gamma}H_{\varepsilon}^{i}\right\Vert _{2}\leq
k\ln\left\vert \varepsilon\right\vert \text{, for\ }\varepsilon
\text{\ small\ enough and }i=1,2.
\]
As $(H_{\varepsilon}^{1}-H_{\varepsilon}^{2})_{\varepsilon}$ is in
${\mathcal{I}}_{L^{2}}(X^{2})$, for a given $m\in\mathbb{N}$, there exists
$C>0$ such that, for $\gamma=\left(  0,0\right)  $, $\gamma=\left(
\alpha,0\right)  $ and $\gamma=\left(  0,\beta\right)  $,
\[
\left\Vert \partial_{x\,y}^{\gamma}(H_{\varepsilon}^{1}-H_{\varepsilon}%
^{2})\right\Vert _{2}\leq C\varepsilon^{m}\text{, for\ }\varepsilon
\text{\ small\ enough.}%
\]
Then, each term in the right hand side of the estimate of $\left\Vert
\partial_{x\,y}^{\left(  \alpha,\beta\right)  }\left(  L_{n,\varepsilon}%
^{1}-L_{n,\varepsilon}^{2}\right)  \right\Vert _{2}$ is less than
$C(k|\ln\varepsilon|)^{n-1}\varepsilon^{m}$, for\ $\varepsilon$%
\ small\ enough. Thus%
\[
\left\Vert \partial_{x\,y}^{\left(  \alpha,\beta\right)  }\left(
L_{n,\varepsilon}^{1}-L_{n,\varepsilon}^{2}\right)  \right\Vert _{2}\leq
Cn(k|\ln\varepsilon|)^{n-1}\varepsilon^{m}\text{, for\ }\varepsilon
\text{\ small\ enough.}%
\]
Finally%
\begin{align*}
\left\Vert \partial_{x\,y}^{\left(  \alpha,\beta\right)  }\left(
S_{\varepsilon}^{1}-S_{\varepsilon}^{2}\right)  \right\Vert _{2} &
\leq C\varepsilon^{m}\sum_{n=1}^{+\infty}\frac{1}{\left(
n-1\right)
!}(k|\ln\varepsilon|)^{n-1}\text{, for\ }\varepsilon\text{\ small\ enough}\\
&  \leq C\varepsilon^{m}e^{k|\ln\varepsilon|}\text{, for\ }\varepsilon
\text{\ small\ enough.}%
\end{align*}
It follows that $\left\Vert \partial_{x\,y}^{\left(  \alpha,\beta\right)
}\left(  S_{\varepsilon}^{1}-S_{\varepsilon}^{2}\right)  \right\Vert
_{2}=O\left(  \varepsilon^{m-k}\right)  $ as $\varepsilon\rightarrow0$ for all
$m\in\mathbb{N}$ and that $(S_{\varepsilon}^{1}-S_{\varepsilon}^{2}%
)_{\varepsilon}$ is in ${\mathcal{I}}_{L^{2}}(X^{2})$. Consequently, $S$ does
not depend on the choice of the representative of $H$ in ${\mathcal{G}}%
_{L^{2}ln}(X^{2})$.

\begin{remark}
The generalization of the results presented in subsections \ref{SbSExpComp}
and \ref{SbSExpGL2} for the exponential, to any entire function is
straightforward by replacing the logarithm scale of growth by an adapted scale.
\end{remark}

\subsection{Properties of the exponential of generalized integral
operators\label{SbSPrExp}}

\begin{proposition}
Let $X$ be an open subset of $\mathbb{R}^{d}$. For $H$ in $\mathcal{G}%
_{C\,ln}(X^{2})$ or in ${\mathcal{G}}_{L^{2}\,ln}(X^{2})$,\ we have%
\[
\widehat{H}\circ e^{\widehat{H}}=e^{\widehat{H}}\circ\widehat{H}%
\ ;\ \ \ \ e^{a\widehat{H}}\circ e^{b\widehat{H}}=e^{(a+b)\widehat{H}%
},\ \text{for all }\left(  a,b\right)  \in\mathbb{R}^{2}\ ;\ \ \ \frac{d}%
{dt}~e^{t\widehat{H}}=\widehat{H}\circ e^{t\widehat{H}}.
\]
Moreover, if $K$ is in $\mathcal{G}_{C\,ln}(X^{2})$ or in ${\mathcal{G}%
}_{L^{2}\,ln}(X^{2})$ and if $\widehat{H}$ and $\widehat{K}$ commute, which
amounts to $%
{\textstyle\int}
H(\cdot_{1},\xi)K(\xi,\cdot_{2})\,\mathrm{d}\xi=%
{\textstyle\int}
K(\cdot_{1},\xi)H(\xi,\cdot_{2})\,\mathrm{d}\xi$, then $e^{\widehat{H}}\circ
e^{\widehat{K}}=e^{\widehat{H}+\widehat{K}}$.
\end{proposition}

\underline{{\bf{Proof}}}. By applying the theorem concerning the
characterization of generalized integral operators by their
kernel, we prove these properties by using the associated kernels.
Denote by $Ker(\cdot)$ the kernel of a generalized integral
operator. From $H$ in ${\mathcal{G}}_{C\,ln}(X^{2})$ or
${\mathcal{G}}_{L^{2}\,ln}(X^{2})$, we define $S$ and $L_{n}$, for
$n\geq1$, as in Theorem-definitions \ref{defexp1} or
\ref{defexp2}. Note that, for all integers $p$, $q$ greater than
$1$, one has
\[
\int_{X}L_{p}(\cdot_{1},\xi)L_{q}(\xi,\cdot_{2})\,\mathrm{d}\xi=L_{p+q}%
(\cdot_{1},\cdot_{2}).
\]
We are going to prove these properties for $H$ in ${\mathcal{G}}_{C\,ln}%
(X^{2})$ and the case where $H$ is in
${\mathcal{G}}_{L^{2}\,ln}(X^{2})$ is treated analogously. We also
don't go back to representatives: All the sums written below are
well defined and independent of representatives. Integrals are
performed on a compact set $\kappa$ such that the support of $H$
is contained in the interior of $\kappa^{2}$\smallskip

$\bullet$ One has $\widehat{H}\circ e^{\widehat{H}%
}=\widehat{H}\circ\widehat{S}+\widehat{H}$ and
\[
Ker\left(  \widehat{H}\circ\widehat{S}\right)  (\cdot_{1},\cdot_{2}%
)=\int_{\kappa}H(\cdot_{1},\xi)S(\xi,\cdot_{2})\,\mathrm{d}\xi=\int_{\kappa
}H(\cdot_{1},\xi)\sum_{n=1}^{+\infty}\frac{L_{n}}{n!}(\xi,\cdot_{2}%
)\,\mathrm{d}\xi=\sum_{n=1}^{+\infty}\frac{L_{n+1}}{n!}(\cdot_{1},\cdot_{2}).
\]
Thus%
\[
Ker\left(  \widehat{H}\circ e^{\widehat{H}}\right)  =\sum_{n=1}^{+\infty}%
\frac{L_{n+1}}{n!}+H.
\]
Knowing that $\int_{\kappa}H(\cdot_{1},\xi)L_{n}(\xi,\cdot_{2})\,\mathrm{d}%
\xi=L_{n+1}(\cdot_{1},\cdot_{2})=\int_{\kappa}L_{n}(\cdot_{1},\xi)H(\xi
,\cdot_{2})\,\mathrm{d}\xi$, one gets
\[
Ker\left(  \widehat{H}\circ e^{\widehat{H}}\right)  =Ker\left(  e^{\widehat
{H}}\circ\widehat{H}\right)  ,
\]
which implies the result.\medskip

$\bullet$ One has $e^{c\widehat{H}}=\widehat{S_{c}}+Id$, for
$c=a,~b,~a+b$, with
\[
S_{c}(\cdot_{1},\cdot_{2})=\sum_{n=1}^{+\infty}\frac{c^{n}}{n!}\int
_{\kappa^{n-1}}H(\cdot_{1},\xi_{1})\cdots H(\xi_{n-1},\cdot_{2})\,\mathrm{d}%
\xi_{1}\cdots\mathrm{d}\xi_{n-1}=\sum_{n=1}^{+\infty}\frac{c^{n}}{n!}%
L_{n}(\cdot_{1},\cdot_{2}).
\]
We have $e^{a\widehat{H}}\circ e^{b\widehat{H}}=\widehat{S_{a}}\circ
\widehat{S_{b}}+\widehat{S_{a}}+\widehat{S_{b}}+Id$, and
\begin{align*}
Ker(\widehat{S_{a}}\circ\widehat{S_{b}})(\cdot_{1},\cdot_{2})  &
=\int_{\kappa}\sum_{n=1}^{+\infty}\frac{a^{n}}{n!}L_{n}(\cdot_{1},\xi
)\sum_{n=1}^{+\infty}\frac{b^{n}}{n!}L_{n}(\xi,\cdot_{2})\,\mathrm{d}\xi\\
&  =\sum_{n=2}^{+\infty}\sum_{k=1}^{n-1}\frac{a^{n-k}b^{k}}{\left(
n-k\right)  !k!}L_{n}(\cdot_{1},\cdot_{2})=\sum_{n=2}^{+\infty}\frac{1}%
{n!}\sum_{k=1}^{n-1}C_{n}^{k}a^{n-k}b^{k}L_{n}(\cdot_{1},\cdot_{2}).
\end{align*}
Thus
\begin{align*}
Ker\left(  \widehat{S_{a}}\circ\widehat{S_{b}}+\widehat{S_{a}}+\widehat{S_{b}%
}\right)   &  =\sum_{n=2}^{+\infty}\frac{1}{n!}\sum_{k=1}^{n-1}C_{n}%
^{k}a^{n-k}b^{k}L_{n}+\sum_{n=1}^{+\infty}\frac{a^{n}}{n!}L_{n}+\sum
_{n=1}^{+\infty}\frac{b^{n}}{n!}L_{n}\\
&  =\sum_{n=1}^{+\infty}\frac{1}{n!}\sum_{k=0}^{n}C_{n}^{k}a^{n-k}b^{k}%
L_{n}=\sum_{n=1}^{+\infty}\frac{1}{n!}\left(  a+b\right)  ^{n}=S_{a+b},
\end{align*}
which implies the required property.\medskip

$\bullet$ One has $e^{t\widehat{H}}=\widehat{S_{t}}+Id$
with $S_{t}=%
{\textstyle\sum\limits_{n=1}^{+\infty}}
\frac{t^{n}}{n!}L_{n}$, so $\frac{d}{dt}e^{t\widehat{H}}=\frac{d}{dt}%
\widehat{S_{t}}=\widehat{\frac{d}{dt}S_{t}}$. Now $\frac{d}{dt}S_{t}=%
{\textstyle\sum\limits_{n=1}^{+\infty}}
\frac{t^{n-1}}{(n-1)!}L_{n}$. Furthermore $\widehat{H}\circ e^{t\widehat{H}%
}=\widehat{H}\circ\widehat{S_{t}}+\widehat{H}$ and
\begin{align*}
Ker\left(  \widehat{H}\circ\widehat{S_{t}}\right)  (\cdot_{1},\cdot_{2})  &
=\int_{\kappa}H(\cdot_{1},\xi)\sum_{n=1}^{+\infty}\frac{t^{n}}{n!}L_{n}%
(\xi,\cdot_{2})\,\mathrm{d}\xi\\
&  =\sum_{n=1}^{+\infty}\frac{t^{n}}{n!}L_{n+1}(\cdot_{1},\cdot_{2})\\
&  =\sum_{n=2}^{+\infty}\frac{t^{n-1}}{(n-1)!}L_{n}(\cdot_{1},\cdot_{2}%
)=\frac{d}{dt}S_{t}(\cdot_{1},y)-L_{1}(\cdot_{1},\cdot_{2}).
\end{align*}
As $L_{1}=H$, one gets $\widehat{H}\circ e^{t\widehat{H}}=\widehat{\frac
{d}{dt}S_{t}}=\frac{d}{dt}e^{t\widehat{H}}$.\medskip

$\bullet$ From $T$ in ${\mathcal{G}}_{C\,ln}(X^{2})$, we define,
for $n\geq1$, $S_{T}$ and $L_{T,n}$ as in Theorem-definition
\ref{defexp1}. Take $H$ and $K$ which commute and $\kappa$ a
compact set such
that the supports of $H$ and $K$ are contained in the interior of $\kappa^{2}%
$. On one hand, we have $e^{\widehat{H}}\circ e^{\widehat{K}}=\widehat{S_{H}%
}\circ\widehat{S_{K}}+\widehat{S_{H}}+\widehat{S_{K}}+Id$ with%
\begin{align*}
Ker\left(  \widehat{S_{H}}\circ\widehat{S_{K}}\right)  (\cdot_{1},\cdot_{2})
&  =\int_{\kappa}\sum_{n=1}^{+\infty}\frac{1}{n!}L_{H,n}(\cdot_{1},\xi
)\sum_{n=1}^{+\infty}\frac{1}{n!}L_{K,n}(\xi,\cdot_{2})\,\mathrm{d}\xi\\
&  =\sum_{n=2}^{+\infty}\sum_{k=1}^{n-1}\frac{1}{\left(  n-k\right)  !k!}%
\int_{\kappa}L_{H,k}(\cdot_{1},\xi)L_{K,n-k}(\xi,\cdot_{2})\,\mathrm{d}\xi\\
&  =\sum_{n=2}^{+\infty}\frac{1}{n!}\sum_{k=1}^{n-1}C_{n}^{k}\int_{\kappa
}L_{H,k}(\cdot_{1},\xi)L_{K,n-k}(\xi,\cdot_{2})\,\mathrm{d}\xi.
\end{align*}
Thus%
\[
Ker\left(  \widehat{S_{H}}\circ\widehat{S_{K}}+\widehat{S_{H}}+\widehat{S_{K}%
}\right)  =\sum_{n=1}^{+\infty}\frac{1}{n!}\sum_{k=1}^{n}C_{n}^{k}\int
_{\kappa}L_{H,k}(\cdot_{1},\xi)L_{K,n-k}(\xi,\cdot_{2})\,\mathrm{d}\xi.
\]
On the other hand, we have $e^{\widehat{H}+\widehat{K}}=Id+\widehat{S_{H+K}}$
with $Ker(\widehat{S_{H+K}})(\cdot_{1},\cdot_{2})=\sum\limits_{n=1}^{+\infty
}\frac{1}{n!}L_{H+K,n}(\cdot_{1},\cdot_{2})$ and, for all $n$ greater than
$1$,
\begin{align*}
L_{H+K,n}(\cdot_{1},\cdot_{2})  &  =\int_{\kappa^{n-1}}\left(  H+K\right)
(\cdot_{1},\xi_{1})\cdots\left(  H+K\right)  (\xi_{n-1},\cdot_{2}%
)\,\mathrm{d}\xi_{1}\cdots\mathrm{d}\xi_{n-1}\\
&  =\sum_{k=1}^{n}C_{n}^{k}\int_{\kappa}L_{H,k}(\cdot_{1},\xi)L_{K,n-k}%
(\xi,\cdot_{2})\,\mathrm{d}\xi.
\end{align*}
This last equality follows from a straightforward induction, which uses mainly
the fact that $H$ and $K$ commute, which implies that $L_{H,p}$ and $L_{K,q}$
have the same property for all integers $p\geq1$ and $q\geq1$. Thus
$Ker(\widehat{S_{H+K}})(\cdot_{1},\cdot_{2})=Ker\left(  \widehat{S_{H}}%
\circ\widehat{S_{K}}+\widehat{S_{H}}+\widehat{S_{K}}\right)  $,
which ends the proof.

\subsection{Example: A unitary generalized integral operator}

In this subsection, we apply the above results to the special case of
operators with symmetrical kernel, which are essential in view of forthcoming
applications to theoretical physics. Fix $X$ an open subset of $\mathbb{R}%
^{d}$.

\begin{propdef}
\label{PropSP}The map $(f,g)\mapsto\int
f(x)\bar{g}(x)\,\mathrm{d}x$ from $\left(
{\mathcal{G}}_{L^{2}}(X)\right)  ^{2}$ to $\overline{\mathbb{C}}$
defines a \emph{generalized scalar product} on
${\mathcal{G}}_{L^{2}}(X)$, that is $(\cdot,\cdot)$ is bilinear,
positive ($(f,f)$ has a representative $\left(
\varphi_{\varepsilon}\right)  _{\varepsilon}$ with $\varphi
_{\varepsilon}\geq0$ for all $\varepsilon\in\left(  0,1\right]  $)
and non degenerate, id est: $(f,f)=0$ in $\overline{\mathbb{C}}$
implies that $f=0$ in ${\mathcal{G}}_{L^{2}}(X)$.
\end{propdef}

\underline{{\bf{Proof}}}. The only non trivial assertion is the
last one. Take $f$ such that $(f,f)=0$ in $\overline{\mathbb{C}}$
and denote by $(f_{\varepsilon})_{\varepsilon}$ one
of its representatives. For all $n\in\mathbb{N}$, $\Vert f_{\varepsilon}%
\Vert_{2}=O(\varepsilon^{n})$ as $\varepsilon\rightarrow0$, that is
$(f_{\varepsilon})_{\varepsilon}$ is in ${\mathcal{I}}_{l^{2}}(X)$.
Furthermore $(f_{\varepsilon})_{\varepsilon}\in{\mathcal{E}}_{L^{2}}(X)$ so
$(f_{\varepsilon})_{\varepsilon}\in{\mathcal{I}}_{l^{2}}(X)\cap{\mathcal{E}%
}_{L^{2}}(X)={\mathcal{I}}_{L^{2}}(X)$ (this result is proved in
\cite{delcroix3} by a method analogous to the one employed in the
proof of
theorem 1.2.3 in \cite{grosser}), which means that $f=0$ in ${\mathcal{G}%
}_{L^{2}}(X)$.


\begin{definition}
We say that a generalized function $H$ of ${\mathcal{G}}_{L^{2}}(X^{2})$ is
\emph{symmetric} if, for all $(x,y)\in X^{2}$, $H(x,y)=\overline{H(y,x)}$ in
$\overline{\mathbb{C}}$.
\end{definition}

\begin{remark}
If $H$ of ${\mathcal{G}}_{L^{2}}(X^{2})$ is symmetric, then $\widehat{H}$ is
symmetric for the generalized scalar product introduced in
proposition-definition \ref{PropSP}, that is $(\widehat{H}\left(  f\right)
,g)=(f,\widehat{H}\left(  g\right)  )$ for all $f,~g\in{\mathcal{G}}_{L^{2}%
}(X^{2})$.
\end{remark}

\begin{definition}
We say that a generalized operator $A$~:~$D(A)\subset{\mathcal{G}}_{L^{2}%
}(X)\rightarrow{\mathcal{G}}_{L^{2}}(X)$ is \emph{unitary} if, for all
$\left(  f,g\right)  \in D(A)^{2}$, one has $(A\left(  f\right)  ,A\left(
g\right)  )=(f,g)$.
\end{definition}

\begin{proposition}
Let $H$ be a symmetric generalized function in ${\mathcal{G}}_{L^{2}%
\,ln}(X^{2})$.\ The generalized integral operator $e^{it\widehat{H}}$ is unitary.
\end{proposition}

\underline{{\bf{Proof}}}. We don't go back to representatives as
in subsection \ref{SbSPrExp}. As $e^{it\widehat{H}}\circ
e^{-it\widehat{H}}=Id$, we have just to prove that
$\left(  e^{it\widehat{H}}(f),g\right)  =\left(  f,e^{-it\widehat{H}%
}(g)\right)  $ in $\overline{\mathbb{C}}$, for all $f,g\in{\mathcal{G}}%
_{L^{2}}(X)$. Let $f$ and $g$ be in ${\mathcal{G}}_{L^{2}}(X)$. By definition
of the exponential, it suffices to show that $\left(  \widehat{S}%
_{it}(f),g\right)  =\left(  f,\widehat{S}_{-it}(g)\right)  $ in $\overline
{\mathbb{C}}$. We have
\begin{align*}
\left(  \widehat{S}_{it}(f),g\right)   &  =\int_{X}\widehat{S}_{it}%
(f)(x)\bar{g}(x)\,\mathrm{d}x\\
&  =\int_{X}\sum_{n=1}^{+\infty}\frac{(it)^{n}}{n!}\int_{X}L_{n}%
(x,y)f(y)\,\mathrm{d}y~\bar{g}(x)\,\mathrm{d}x\\
&  =\sum_{n=1}^{+\infty}\frac{(it)^{n}}{n!}\int_{X}\int_{X}L_{n}%
(x,y)f(y)\bar{g}(x)\,\mathrm{d}y\mathrm{d}x.
\end{align*}
As $H$ is symmetrical, so is $L_{n}$ for all $n\geq1$ and, by applying
Fubini's theorem, one gets
\begin{align*}
\left(  \widehat{S}_{it}(f),g\right)   &  =\sum_{n=1}^{+\infty}\frac{(it)^{n}%
}{n!}\int_{X}\int_{X}\bar{L}_{n}(y,x)f(y)\bar{g}(x)\,\mathrm{d}x\mathrm{d}y\\
&  =\int_{X}f(y)\left(  \int_{X}\sum_{n=1}^{+\infty}\frac{(it)^{n}}{n!}\bar
{L}_{n}(y,x)\bar{g}(x)\,\mathrm{d}x\right)  \,\mathrm{d}y\\
&  =\left(  f,\widehat{S}_{-it}(g)\right)  .
\end{align*}

\end{document}